\documentclass[final]{siamltex}
\usepackage{amsfonts}
\usepackage{epsf,graphicx,epsfig}

 \newtheorem {remark}{Remark}
 \newtheorem {assumption}{Assumption}

\def \Z {{\mathbb{Z}}}

\def \N {{\mathbb{N}}}
\def \R {{\mathbb {R}}}
\def \expec{\mathbb{E}}
\def \T{\mathbb{T}}

\def\bra{\langle}
\def\cet{\rangle}

\def \bov {{\bf {v}}}
\def \boV {{\bf {V}}}
\def \bou {{\bf {u}}}
\def \boU {{\bf {U}}}
\def \bom {{\bf {m}}}
\def \boM {{\bf {M}}}
\def \boA {{\bf {A}}}
\def \boE {{\bf {E}}}

\def \onef{{\bf 1}}

\newcommand{\norm}[1]{\left\|#1 \right\|}

\newcommand{\abs}[1]{\left| #1 \right|}

\newcommand{\ep}{\epsilon}

\newcommand{\ran}{{\rm Ran}}

\newcommand{\Pro}{\mathbb{P}}
\newcommand{\Ec}{{\mathcal E}}
\newcommand{\Ll}{{\mathcal L}}

\newcommand{\td}{D_q}

\newcommand{\Fe}{F_{\ep}}
\newcommand{\Fbe}{F^\alpha_\ep}
\newcommand{\Fda}{F_{\ep,n}^\alpha}

\newcommand{\Fat}{F_{\ep}^{AT}}

\newcommand{\MS}{MS}

\newcommand{\at}{AT_\ep}
\newcommand{\pone}{L_\ep}
\newcommand{\ptwo}{S^q_\ep}

\newcommand{\gsbv}{GSBV}
\newcommand{\sbv}{SBV}
\newcommand{\bv}{BV}

\newcommand{\weakly}{\rightharpoonup}
\newcommand{\aplim}{{\rm aplim}}
\newcommand{\Glim}{\Gamma\!\!-\!{\rm lim}}

\newcommand{\eqref}[1]{(\ref{#1})}

\title{Hierarchical models in statistical inverse problems and the Mumford--Shah functional
\thanks{The first author was supported by Emil Aaltonen foundation, Graduate School of 
Inverse Problems (Academy of Finland) and Finnish
Centre of Excellence in Inverse Problems Research (Academy of Finland).
Authors are thankful for Luca Rondi and Samuli Siltanen for helpful comments.}}

\author{
Tapio Helin\thanks{Department of Mathematics and System Analysis, Helsinki University of Technology, P.O.Box 1100, FI-02015 TKK, Finland ({\tt tapio.helin@tkk.fi})}.
\and Matti Lassas\thanks{Department of Mathematics and Statistics,
University of Helsinki, P.O. Box 68 (Gustaf Hallstromin katu 2b) FI-00014, Finland ({\tt matti.lassas@helsinki.fi}.)}}

\begin{document}

\maketitle

\begin{abstract}
The Bayesian methods for linear inverse problems is studied using hierarchical Gaussian models.
The problems are considered with different discretizations, and we analyze
the phenomena which appear
when the discretization becomes finer. A  hierarchical
solution method for signal restoration problems
is introduced and studied
with arbitrarily fine discretization. 
We show that the maximum a posteriori estimate converges 
to a minimizer of the Mumford--Shah
functional, up to a subsequence. 
A new result regarding the existence of a minimizer of 
the Mumford--Shah functional is proved.

Moreover, we study the inverse problem under different assumptions on the 
asymptotic behavior of the noise as discretization becomes finer. 
We show that the maximum a posteriori and conditional mean estimates
converge under different conditions.
\end{abstract}



\pagestyle{myheadings}
\thispagestyle{plain}




\section{Introduction}


We study hierarchical Bayesian methods for linear inverse problems.
In particular, we consider inverse problems with different discretizations and 
the phenomena which appear when the
discretization is refined. The effect of fine discretization has recently been studied
for Gaussian inverse problems in \cite{KSpaper,Lasanen,NP}, 
and motivated by this development we consider hierarchical Gaussian
models. More precisely, we introduce a hierarchical 
solution method and analyze its properties with arbitrarily fine
discretization. 

The inverse problem we consider is the  linear signal restoration problem where the measurement $m(t)$ 
relates indirectly to the unknown signal $u(t)$ via 
\begin{equation}
  \label{eq:linear_ip}
  m(t) = Au(t) + e(t),\quad t\in \T.
\end{equation}
Here, $\T$ is the unit circle which we frequently consider as the interval $[0,1]$ with the end points identified. Furthermore, $A$ is a smoothing linear integral operator, and $e(t)$ is random noise. 
The signals are considered on the unit circle $\T$ to avoid the complicated boundary effects 
that fall outside the scope of this paper.

In the Bayesian approach $u(t)$ and $e(t)$ are modelled as random functions.
Let us denote by $U(t,\omega)$ and $\Ec(t,\omega)$ random functions where $\omega\in\Omega$
is an element of a complete probability space $(\Omega,\Sigma,\Pro)$ and $t\in\T$. 
The distribution of $U(t,\omega)$ and $\Ec(t,\omega)$ model our {\em a priori} knowledge on the
unknown signal  $u(t)$ and error $e(t)$, respectively, before the measurement is obtained.
Below, the variable
$\omega$ is often omitted. 
The ideal measurement is considered to be a realization of the random function $M(t) = AU(t)+\Ec(t)$ on
$t\in \T$.
In Bayesian inversion the aim
is to make statistical inference on $U$ given a realization $m$ of the
random function $ M$, and the Bayesian solution  to the inverse problem means  finding
the conditional probability distribution of $U$,
called the {\em posterior distribution}, or some estimates for this distribution. 
Typically studied point estimates are the expectation of the 
posteriori distribution called the {\em conditional mean} (CM) estimate and
the maximum point of the posterior density called the {\em maximum a posteriori} (MAP) estimate.



In Bayesian inversion a reconstruction method is said to be edge-preserving 
if the functions $u$ which have high
probability with respect to the posterior distribution
are roughly speaking piecewise smooth and have rapidly changing values only in a set of small measure.
In the finite-dimensional Bayesian inversion theory a number of methods have been introduced for obtaining
edge-preserving reconstructions \cite{boumansauer,chantas,GG,raj}. 
In this paper 
the prior distribution of the random function $U$ has a Gaussian distribution such that
its covariance depends on an auxiliary random
function $V$. Moreover, the random function $V$ has a Gaussian distribution.
Such a model is called {\em hierarchical Gaussian model}.
With a fixed discretization similar models have been studied in inverse problems in \cite{VKK}.
Furthermore, in the work by Calvetti and Somersalo
\cite{CS1,CS2} hierarchical methods have been used for image processing problems to obtain 
edge-preserving and numerically efficient reconstruction algorithms.
We also mention that the edge-preserving reconstruction methods have been extensively studied
in the deterministic problem setting, see e.g., \cite{candesdonoho,CBAB,mumfordshah,PM,ROF,Sethian,chan}.
Our main result in this paper connects computing the
MAP estimate of a hierarchical Gaussian model to
the minimization of the Mumford--Shah functional \cite{mumfordshah} 
used in image processing. As a byproduct we also present new results concerning the existence of a minimizer
of the Mumford--Shah functional.

Let us next discuss the discretization of Bayesian inverse problems. Above we have
considered $U(t)$ and $M(t)$ as random functions defined on the unit circle. 
For any practical computations such models have to be discretized, i.e., to be approximated
by random variables taking values in a finite dimensional space. Roughly speaking, a Bayesian
model is said to be {\em discretization invariant}, if for fixed model parameters it works coherently
at any level of discretization.
For an extended discussion on the discretization invariance and 
the relation of the practical measurement models and the computational models considered below, see \cite{LSS}. 

In the ideal model the noise $\Ec$ can be considered as a background
noise. In this paper we will further assume that the practical measurement setting
produces an additional instrumentation noise.
More precisely, we assume that the practical measurement 
can be modelled as a 
realizations of a random variable
$M_k = P_k M+E_k$, $k=1,2,3...$,
where operator $P_k$ is a finite dimensional projection. The random variable
$E_k$ describes the instrumentation noise and it takes
values in the range of $P_k$.
Increasing the number $k$ corresponds here to the case when we make more or finer observations
of the ideal measurement signal $M(t)$.
Moreover, in practical computations also $U$ needs to be approximated by a finite dimensional random variable $U_n$ which leads us to consider the {\em computational model}
\begin{equation}
  \label{eq:compmodel}
  M_{kn} = A_k U_n +\Ec_k
\end{equation}
where $k,n\in\N$ are parameters related to discretizing the measurement and the unknown, respectively.
In equation \eqref{eq:compmodel} we have $A_k = P_k A$ and $\Ec_k$ is a random variable
in range of $P_k$ satisfying
\begin{equation}
  \Ec_k = P_k \Ec + E_k.
\end{equation}


In developing new Bayesian algorithms, it is important to study if the posterior distribution
given by problem \eqref{eq:compmodel} or some preferred estimate converges when $k$ or $n$
increase.
This question is often non-trivial. For example, for the total variation prior it is proven in \cite{LS} 
that the MAP and CM estimates converge under different conditions as discretization is refined. Moreover,
if the free parameters of the discrete total variation priors are chosen so that the posterior
distribution converges, then the limit is a Gaussian distribution.
Hence, the key property of the total variation prior is lost in very fine discretizations.
This example illustrates the difficulty involved in discretizing non-Gaussian distributions.
Also in this paper we will observe that the convergence of the MAP and CM estimates 
occurs in different cases.

 
Let us next formally define the discrete models we study.
Set $N=2^n$ and let points $t_j = j/N$, $j=0,1,...,N$, and $t_0$ identified with $t_N$, denote an
equispaced mesh on $\T$.
We define $PL(n)$ to be the space of continuous functions $f\in C(\T)$ such that $f$
is linear on each interval $[t_j,t_{j+1}]$ for $0\leq j<N$.
Furthermore, let $PC(n)$ be the space of functions $f\in L^2(\T)$ such that $f$ is constant
on each interval $[t_j,t_{j+1})$ for $0\leq j<N$.
Denote by $Q_n : L^2(\T)\to PC(n)$ the orthogonal projections with respect to 
$L^2(\T)$ inner product and let $Q = Q_0$ be the projection to constant functions.
Define the operator $D_q=D+\ep^qQ : H^1(\T)\to L^2(\T)$ where 
$\epsilon>0$, $q>1$ and $D$
is the derivative with respect to $t\in\T$ .

The hierarchical structure is defined in two steps. First, let $V_{n,\ep}$ be a Gaussian
random variable in $PL(n)$ with density function
\begin{equation}
  \label{eq:formalprior1}
  \Pi_{V_{n,\epsilon}}(v)=
  c \exp\left( -\frac{N^{\alpha}}2 \int_\T  \left(\frac{1}{4\ep}|v(t)-1|^2+\ep |Dv(t)|^2\right) \,dt\right)
\end{equation}
where $v\in PL(n)$, $\alpha\in\R$ and $N$ is the number of mesh points. 
Here and below $c$ is a generic constant whose value may vary.
Then choose $v_{n,\ep}$ to be a sample of $V_{n,\ep}$. The random variable $U_{n,\ep}$, conditioned on $v_{n,\ep}$,
is then defined as a Gaussian random variable on $PL(n)$ with density function
\begin{equation}
  \label{eq:formalprior2}
  \Pi_{U_{n,\ep} \; | \; V_{n,\ep}}(u \; | \; v_{n,\ep})=c' \exp\left( -\frac{N^{\alpha}}2 \int_\T
  (\epsilon^2+|Q_nv_{n,\ep}(t)|^2) {|D_q u(t)|^2}\,dt\right)
\end{equation}
where $u\in PL(n)$. Note that the constant $c'$ depends on $v_{n,\ep}$. 
Since the density function presentation in finite-dimensional Hilbert spaces is non-standard,
we give in Section 2 a definition 
of random variables $U_{n,\ep}$ and $V_{n,\ep}$ based on the coordinate representation.

Roughly speaking, the sample $v_{n,\epsilon}$ has a high probability if it varies from $\onef$ only little
and this variation becomes less smooth if $\ep$ is decreased.
A sample of $U_{n,\epsilon}$ has a high probability if it
varies rapidly only near the points where $v_{n,\epsilon}$ is close to zero. 
Hence the role of parameter $\ep>0$ is to control how sharp jumps $U_{n,\epsilon}$ can have and
consequently, we call it the {\em sharpness} of the prior.
Furthermore, the parameter $\alpha$ describes the {\em scaling} of the prior information.
The bigger the value of $\alpha$ the more concentrated the prior distribution.

In consequence of the construction above the probability density of the joint distribution of $(U_{n,\epsilon},V_{n,\epsilon})$ 
has a form
\[
  \Pi_{(U_{n,\ep},V_{n,\ep})}(u,v)
  =c \exp\left(-\frac{N^{\alpha}}2 F^\alpha_{\ep,n}(u,v)\right)
\]
where $(u,v)\in PL(n)\times PL(n)$ and
\begin{eqnarray}
  \label{eq:fn_formal}
  F^\alpha_{\ep,n}(u,v) & = &
  \int_\T \bigg(-N^{1-\alpha} \log(\epsilon^2+|Q_nv|^2) + \nonumber \\
  & & \quad +(\epsilon^2+|Q_nv|^2)|\td u|^2 + \ep |D v|^2 + \frac{1}{4\ep} |1-v|^2 \bigg) dt.
\end{eqnarray}
The logarithmic term in \eqref{eq:fn_formal} appears due to the fact that the 
normalization constant $c'$ in \eqref{eq:formalprior2} depends on $v_{n,\epsilon}$.
It turns out that the functional $F^\alpha_{\ep,n}$ is closely connected to the widely studied 
segmentation method in deterministic image and signal processing,
namely, the Mumford-Shah functional
\begin{equation}
  \label{eq:strong_ms}
  F(u,K) = \int_{\T\setminus K} |Du|^2 dt + \sharp(K) + \int_\T |Au-m|^2 dt
\end{equation}
with respect to function $u$ and the set $K$ of the points where $u$ jumps \cite{mumfordshah}.
The notation $\sharp(K)$ stands for the number of points in $K$. 
This functional is known to be difficult to handle numerically and
a number of approximations to the variational problem of minimizing \eqref{eq:strong_ms} have been introduced.
In \cite{AT1,AT2} it is shown that the Mumford--Shah functional 
can be approximated by elliptic functionals
in the sense of  $\Gamma$-convergence.
These Ambrosio-Tortorelli functionals are the key element in our presentation.


Let us describe our main results.
We study the behaviour of the MAP estimate 
in the case when the discretization parameters $k$ and $n$ are coupled.
For the sake of presentation we assume $k=n$ and drop $k$ from the notations.
Furthermore, we assume that $\Ec_n$ is white noise with variance depending on $n$ and
scaling parameter $\kappa$.
More precisely, $\Ec_n$ is
a Gaussian random function on $\T$ taking values in $\ran(P_n)$ with zero expectation and covariance
\begin{equation}
	\label{intro_noise}
  \expec \left(\bra \Ec_n,\phi\cet_{L^2} \bra \Ec_n,\psi\cet_{L^2}\right) = N^{-\kappa}\bra \phi,\psi\cet_{L^2}
\end{equation}
for any $\phi,\psi\in \ran(P_n)$. Notice that consequently $\expec \norm{\Ec_n}^2_{L^2} = N^{1-\kappa}$
and the choice of $\kappa$ describes how the norm of the noise is expected to behave asymptotically.
We emphasize that the case when $\kappa > 1$ corresponds an assumption that more
measurements produces {\em better accuracy expectation} whereas with $\kappa=1$ one assumes that
the accuracy in the norm of $L^2(\T)$ is expected to stay stable.
An example of the case $\kappa\geq 1$ is when the background noise
$\Ec$ is negligible and the instrumentation noise follows asymptotics \eqref{intro_noise}.
The case $\kappa=0$ corresponds to the discretization of the Gaussian white noise,
see \cite{LSS}.
To be able to prove positive results for the convergence of the MAP estimates we will assume
\begin{equation}
	\label{k_equals_a}
  \kappa = \alpha.
\end{equation}
This implies that the scaling parameter of the prior distribution is determined by the variance of the noise
in discretized measurements. The case when \eqref{k_equals_a} is not valid is discussed in Remark \ref{k_not_equal_a}.
Due to the equality \eqref{k_equals_a} we drop the notation $\kappa$ and use $\alpha$ as the scaling parameter of
the noise distribution.

Under these assumptions the MAP estimate for  $(U_{n,\epsilon},V_{n,\epsilon})$ corresponding the measurement $m_n$, $\lim_{n\to\infty} m_n = m$ in $L^2(\T)$, is a minimizer 
\begin{equation}
  \label{eq:maps}
   \left(u_{n,\epsilon}^{MAP},v_{n,\epsilon}^{MAP}\right) 
  \in {\rm argmin}_{(u,v)\in PL(n)\times PL(n)} 
  \left(F_{\ep,n}^\alpha(u,v) + \norm{A_n u - m_n}^2_{L^2} \right).
\end{equation}
In the Theorems \ref{neg_result} and \ref{pos_result} we prove for the MAP estimates:
\begin{itemize}
\item[(a)] For $\alpha=0$ the minimization problems \eqref{eq:maps} diverge as $n\to\infty$.
\item[(b)] For $\alpha\geq 1$ the MAP estimates $ \left(u_{n,\epsilon}^{MAP},v_{n,\epsilon}^{MAP}\right) $
converge to the minimizer, denoted $\left(u_{\epsilon}^{MAP},v_{\epsilon}^{MAP}\right)$, of 
a perturbed Ambrosio--Tortorelli functional as $n\to \infty$. 
Moreover, the functions $\left(u_{\epsilon}^{MAP},v_{\epsilon}^{MAP}\right)$ are shown to 
converge up to a subsequence to the minimizer of the Mumford--Shah functional  (\ref{eq:strong_ms}) as $\ep\to 0$.
\end{itemize}

\noindent In \cite{tapio} and Remark \ref{remarkb} 
the following is shown, with slightly different assumptions on operator $A$, for the convergence of prior distributions and the CM estimate
\begin{itemize}
\item[(a')] For $\alpha=0$ the random variables $(U_{n,\epsilon},V_{n,\epsilon})$ converge in distribution
on $L^2(\T)\times L^2(\T)$ and the CM estimates $(U^{CM}_{n,\epsilon},V^{CM}_{n,\epsilon})$
converge in $L^2(\T)\times L^2(\T)$ as $n\to\infty$.
\item[(b')] For $\alpha\geq 1$ the random variables $(U_{n,\epsilon},V_{n,\epsilon})$
converge to zero as $n\to\infty$.
\end{itemize}
The type of convergence in (b') is discussed in Remark \ref{remarkb}.
Consequently,
the results (a),(b) and (a'),(b') illustrate how the convergence properties of the MAP and CM estimates
are different for hierarchical Gaussian models.

Let us recall that the CM and MAP estimate coincide for finite dimensional Gaussian inverse problems \cite{KS}. 
Typically the MAP estimates are computationally faster to obtain than the CM estimates and
thus in inverse problems close to Gaussian ones
the MAP estimate is used as an approximation of the CM estimate.
The above results show that this is not the case for the hierarchical Gaussian models in general.


Finally, let us consider the current perspectives to Bayesian modelling and how this paper connects to earlier work.
Bayesian inversion in infinite-dimensional function spaces were first studied by Franklin in \cite{Franklin}.
This research has then been continued and generalized in \cite{Fitzpatrick, LPS, Mandelbaum}.
The convergence of the posterior distribution is studied in \cite{tapio, Lasanen, LSS, Piiroinen}.
In relation to result $(b)$ the convergence
of posterior distribution is studied in \cite{HP,HP2,NP}
when objective information becomes more accurate with Gaussian prior and noise distributions.
For a general resource on the Bayesian inverse problems theory and computation see \cite{CS3, KS}.
For non-Gaussian noise models in statistical inverse problems see \cite{Hohage}.
The Mumford--Shah functional has been applied to inverse problems for instance in \cite{ramlauring,rondi, rondisantosa} and for related work in image processing problems
see \cite{AV,candestao,chambolle,chanshen}.
Finally we mention that variational approximation with $\Gamma$-convergence is used earlier 
in the context of inverse problems in e.g. \cite{GKLU,LS,rondi,rondisantosa}. 

This paper is organized as follows. In Section 2 we introduce the stochastic model and
necessary tools to tackle the convergence problems related to MAP estimates.
Section 3 covers the main results and the proofs are postponed to Sections 5 and 6.
In Section 4 we study the existence of MAP estimates and
Section 5 discusses the cases when desired convergence does not take place.
In Section 6 the proofs related $\Gamma$-convergence and equi-coerciveness of the functionals.
Finally, in Section 7 we illustrate the method in practise by numerical examples.


\section{Definitions}


In this section we cover the stochastic model introduced in \cite{tapio} and furthermore give
the main tools and theoretical results concerning the variational problem of the MAP estimate.
Let us first introduce some notation. 
Most function spaces in our presentation have structure of a real separable Hilbert space.
We often use the $L^2$-based Sobolev spaces $H^s(\T)$ for any $s\in\R$ equipped with Hilbert space
inner product
\[
  \bra \phi,\psi \cet_{H^s} := \int_\T ((I-\Delta)^{s/2}\phi)(t) ((I-\Delta)^{s/2} \psi)(t) dt
\]
for any $\phi,\psi\in H^s(\T)$ where $\Delta = \frac{d^2}{dt^2}$. 
However, we also study the Banach structure of $H^s(\T)$ with
dual space $H^{-s}(\T)$. In this setting the Banach dual pairing is denoted by 
$\bra \cdot, \cdot \cet_{H^{-s}\times H^s}$. We also denote $H^s(\T; [a,b])$, $s\geq 0$,
for functions $f\in H^s(\T)$ such that $a\leq f \leq b$ a.e. for $a,b\in\R$.
Furthermore, we discuss spaces $H^s(a,b) = \{f \; | \;  f = g|_{(a,b)}, g\in H^s(\R)\}$ for $a,b,s\in\R$.
We say that a sequence $\{x_j\}_{j=1}^\infty$ converges to $x$ strongly in Banach space $X$,
if $\lim_{j\to\infty} \norm{x-x_j}_X=0$ as $j\to\infty$.

Recall from Hilbert space valued stochastics \cite{boga} that a covariance operator $C_X$ 
of a Gaussian random variable $X:\Omega\to H$ is defined by equality
\[
  \expec \left(\bra X-\expec X,\phi\cet_H \bra X-\expec X,\psi\cet_H\right) = \bra C_X \phi,\psi\cet_H
\]
for all $\phi,\psi\in H$. We call a Gaussian random variable centered if $\expec X = {\bf 0}$.


We use a perturbed derivative $\td = D+\ep^q Q$,
where $D= \frac d{dt}$ and
\[
	(Qf)(t)=\left(\int_\T f(t')dt'\right)\onef(t)
\]
for $\onef(t) \equiv 1$ and any $f\in L^1(\T)$. This construction guarantees that
$\td : H^1(\T)\to L^2(\T)$ and $\td |_{PL(n)} : PL(n) \to PC(n)$ are invertible mappings.


\subsection{Bayes modelling}
\label{sec:bayes}

Let us now shortly describe how we define the Bayesian maximum a posteriori estimate for
the computational model given in equation \eqref{eq:compmodel}.
Let $(H_1,\bra\cdot,\cdot\cet_1)$ and $(H_2,\bra\cdot,\cdot\cet_2)$
be two real Hilbert spaces such that $\dim H_1 = J$ and $\dim H_2 = K$.
Assume that $U_n$ obtains realizations in a $H_1$
and the range of the measurement projection $P_k$ is $H_2$.
Furthermore, let ${\mathcal I} : H_1\to \R^J$ and ${\mathcal J} : H_2 \to \R^K$
be two arbitrary isometries.
Let us now map equation \eqref{eq:compmodel} to a matrix equation
\begin{equation}
  \boM_{kn} = {\mathcal J} M_{kn} = \boA_{kn} \boU_n + \boE_k
\end{equation}
where $\boA_{kn} = {\mathcal J} A_{k} {\mathcal I}^{-1} \in \R^{K\times J}$, $A_k = P_kA$,
$\boE_k = {\mathcal J} E_k$ and 
$\boU_n = {\mathcal I}U_n : \Omega \to \R^J$.
If the a priori and likelihood distributions above 
are absolutely continuous with respect to Lebesgue measure
the posteriori distribution can be obtained by the Bayes formula:
the posteriori density $\pi_{kn}$ then has the form
\begin{equation}
  \label{eq:bayesformula}
  \pi_{kn} (\bou \; | \; \bom) = \frac{\Pi_n(\bou)\Gamma_{kn}(\bom \; | \; \bou)}{\Upsilon_{kn}(\bom)}
\end{equation}
where $\bou \in \R^J$ and $\bom \in \R^K$.
In equation \eqref{eq:bayesformula} functions
$\Pi_n$ and $\Gamma_{kn}$ are the prior and likelihood densities, respectively, and
$\Upsilon_n$ is the density of $\boM_{kn}$ \cite{KS}. 
The standard definition of the maximum a posteriori
estimate is then 
\[
  \bou^{MAP}_{kn} \in {\rm argmax}_{\bou \in \R^n} \pi_{kn}(\bou \; | \; \bom)
\]
where the set on the right-hand side consist of all points $\bou$ maximizing $\pi_{kn}(\cdot \; | \; \bom)$.
The value of
\[
  u^{MAP}_{kn} = {\mathcal I}^{-1}\left(\bou^{MAP}_{kn}\right) \in H_1
\]
is commonly defined as the MAP estimate of problem \eqref{eq:compmodel}.
Another point estimate in Bayesian inversion is the CM estimate
which is defined as the integral
\[
  \bou^{CM}_{kn} = \int_{\R^N} \bou \pi_{kn}(\bou \; | \; \bom) d\bou \quad {\rm and} \quad
  u^{CM}_{kn} = {\mathcal I}^{-1} \left(\bou^{CM}_{kn}\right)\in H_1.
\]
We note that although the posterior density depends on the inner products $\bra\cdot,\cdot\cet_1$
and $\bra \cdot,\cdot\cet_2$ both point estimates are invariant with respect to such choices.
For more information about the point estimates see \cite{LSS} for CM
and \cite{Hegland} for MAP estimation in Hilbert spaces.


\subsection{The prior model}
\label{sec:prior_model}

In this subsection we introduce the prior model discussed in the introduction
and explain the meaning of the density function representation in equations \eqref{eq:formalprior1} and
\eqref{eq:formalprior2}.
Let us first review the infinite dimensional prior model introduced in \cite{tapio}.
Consider a centered Gaussian distribution $\lambda^v$ on $L^2(\T)$ with covariance
operator
\[
  C_U(v) = L \Lambda(v) L^* : L^2(\T) \to L^2(\T)
\]
with $L = \td^{-1}:L^2(\T)\to L^2(\T)$ and multiplication operator $\Lambda(v):L^2(\T)\to L^2(\T)$ defined as
\[
  (\Lambda(v)f)(t) = \frac{1}{\ep^2+v(t)^2}\cdot f(t), \quad t\in\T,
\]
for any $v\in L^2(\T)$. Let us now formally discuss the qualitative behavior of $\lambda^v$.
Such a distribution has the property that in a set of $t\in\T$ where $v(t)^2$ is large
the samples from distribution $\lambda^v$ are likely to be smooth. Vice versa, in sets where $v(t)^2$ is 
small the distribution allows more rapid changes. 

Next we set the prior distribution of random variable $U$ to be $\lambda^v$.
However, the crucial step in hierarchical modelling is to model the values of $v$ with a random variable $V$.
Thus, instead of knowing the exact locations of the 
jumps, we model how they are distributed. 
In \cite{tapio} the random variable $V$ 
is Gaussian with expectation $\expec V = \onef$ and
covariance operator $C_V = \left(\frac 1{4\ep} I -\ep \Delta\right)^{-1}$
on $L^2(\T)$. Denote the distribution of $V$ on $L^2(\T)$ by $\nu$. 
The joint distribution of the random variable $(U,V) : \Omega\to L^2(\T)\times L^2(\T)$ is then defined to satisfy
\begin{equation}
  \label{infdimmodel}
  \lambda(E\times F) = \int_{F} \lambda_v(E) d\nu(v)
\end{equation}
for any Borel measurable sets $E,F\subset L^2(\T)$. This construction is shown in \cite{tapio} to be well-defined.

In the following we define the finite dimensional prior structure studied in this paper with all scalings $\alpha\in \R$.
In \cite{tapio} these random variables are shown to converge to $U$ and $V$ in distribution on 
$L^2(\T) \times L^2(\T)$ when $\alpha=0$.
First define two inner products on $H^1(\T)$, namely,
\begin{equation}
  \label{eq:innerprods}
  \bra f,g \cet_1 := \bra \td f, \td g\cet_{L^2} \quad {\rm and} \quad
  \bra f,g \cet_2 := \bra C_V^{-\frac 12} f,C_V^{-\frac 12} g\cet_{L^2}
\end{equation}
for any $f,g\in H^1(\T)$. Next construct two orthonormal basis 
$\{f_j\}_{j=1}^\infty,\{g_j\}_{j=1}^\infty\subset H^1(\T)$
with respect to inner products $\bra\cdot,\cdot\cet_1$ and 
$\bra\cdot,\cdot\cet_2$, respectively, in the following way:
for any $n\in\N$ we have
$\{f_j\}_{j=1}^N, \{g_j\}_{j=1}^N \subset PL(n)$
where $N = 2^n$, $n\in \Z_+$. 
Such a construction can be obtained, e.g., using the Gram-Schmidt orthonormalization procedure.
To simplify our notations we assume that the probability space $(\Omega,\Sigma,\Pro)$ has the additional structure
$\Omega = \Omega_1\times \Omega_2$, $\Sigma = \overline{\Sigma_1 \otimes \Sigma_2}$ and $\Pro = \Pro_1\otimes \Pro_2$. We denote $\omega = (\omega_1,\omega_2)\in \Omega_1\times \Omega_2$.

\begin{definition}
  Define $V^\alpha_{n,\ep} : \Omega_2 \to PL(n)\subset L^2(\T)$ as
  \[
    V^\alpha_{n,\ep}(\omega_2) = \sum_{j=1}^N \boV^\alpha_{j,N,\ep} (\omega_2) g_j + \onef
  \]
  where the random vector $\boV_{n,\ep}^\alpha = (\boV_{j,N,\ep}^\alpha)_{j=1}^N : \Omega_2\to \R^N$ are a centered Gaussian random variable
  with covariance matrix $C_{\boV_{n,\ep}^\alpha} = N^{-\alpha} I\in \R^{N\times N}$.
\end{definition}

\begin{definition}
  Let $U_{n,\ep}^\alpha : \Omega \to L^2(\T)$ be the random variable
  \[
    U^\alpha_{n,\ep}(\omega_1,\omega_2) = \sum_{j=1}^N \boU^\alpha_{j,N,\ep}(\omega_1,\omega_2) f_j
  \]
  where the random vector $\boU_{n,\ep}^\alpha(\omega) = (\boU^\alpha_{j,N,\ep})_{j=1}^N \in \R^N$ is given the following
  structure: Denote by $\omega_2 \mapsto {\bf C}(\omega_2) \in \R^{N\times N}$ a random matrix such that
  \[
    {\bf C}_{jk}(\omega_2) = N^{-\alpha}\bra \Lambda_n(V^\alpha_{n,\ep}(\omega_2)) \td f_j,\td f_k\cet_{L^2}.
  \]
  Due to the positive definiteness of ${\bf C}$ we can define
  \[
    \boU^\alpha_{n,\ep}(\omega) = {\bf C}(\omega_2)^{\frac 12} {\bf W}^N(\omega_1)
  \]
  where ${\bf W}^N : \Omega_1\to \R^N$ is a centered Gaussian random variable with identity covariance matrix.
\end{definition}

Following the procedure shown in Section \ref{sec:bayes} 
choose ${\mathcal I}_1,{\mathcal I}_2 : PL(n)\to \R^N$
to be two isometries
with respect to inner products $\bra\cdot,\cdot\cet_1$ and $\bra \cdot,\cdot\cet_2$, respectively, and with
the usual inner product of $\R^N$.
Clearly, it follows that the vector presentation is then
\[
  \boU^\alpha_{n,\ep} = {\mathcal I}_1 U^\alpha_{n,\ep} \quad {\rm and} \quad
  \boV^\alpha_{n,\ep} = {\mathcal I}_2 (V^\alpha_{n,\ep}-\onef).
\]
In \cite{tapio} it was shown that if
$u, v\in PL(n)$ are arbitrary and $\bou={\mathcal I}_n u, \bov = {\mathcal J}_n v \in \R^N$ 
then it holds that the probability density function of $\boV_{n,\ep}^\alpha$ in $\R^N$ is
\begin{equation}
  \label{eq:vdens}
  \Pi_{\boV^\alpha_{n,\ep}}(\bov) = c
  \exp\left(-\frac{N^{\alpha}}{2} \left(\ep\norm{Dv}^2_{L^2}+
    \frac{1}{4\ep}\norm{v-\onef}^2_{L^2}\right)\right)
\end{equation}  
and the conditional probability density function of $\boU_{n,\ep}^\alpha$ in $\R^N$ is
\begin{eqnarray}
  \label{eq:udens}
   \Pi_{\boU^\alpha_{n,\ep} | \boV^\alpha_{n,\ep}} (& \bou &\; | \; \bov) = \nonumber \\
  \qquad \qquad  & c & \exp \left(-\frac{N^{\alpha}}{2}  \left(\int_\T \left(-N^{1-\alpha}\log(\ep^2+(Q_n v)^2) + (\ep^2+(Q_n v)^2)|\td u|^2\right) dt\right)\right).
\end{eqnarray}
With the same assumptions the joint prior density then takes the form 
\[
  \Pi_{(\boU^\alpha_{n,\ep},\boV^\alpha_{n,\ep})}(\bou,\bov) 
  = \Pi_{\boU^\alpha_{n,\ep} | \boV^\alpha_{n,\ep}}(\bou\; | \; \bov) \cdot \Pi_{\boV^\alpha_{n,\ep}}(\bov)
  = c \exp\left(-\frac{N^{\alpha}}2 \Fda(u,v)\right)
\]
where the functional $\Fda$ is given in the following definition.
\begin{definition}
  For any $\ep>0$, $n\in\N$, and $\alpha\in\R$ let $\Fda : H^1(\T) \times H^1(\T) \to \R \cup \{\infty\}$ be functional such that
  \begin{eqnarray*}
    \Fda(u,v) & = & \int_\T \bigg(-N^{1-\alpha}\log(\ep^2+(Q_nv)^2) \\
   & & \quad + (\ep^2+(Q_nv)^2)|\td u|^2+\ep |Dv|^2 + \frac{1}{4\ep} (1-v)^2\bigg) dx
  \end{eqnarray*}
  when $(u,v)\in PL(n)\times PL(n)$ and $\Fda(u,v) = \infty$ when $(u,v)\in (H^1(\T)\times H^1(\T)) \setminus (PL(n)\times PL(n))$.
\end{definition}


\subsection{Variational approximation and the functions of bounded variation}
\label{sec:Glim_and_BV}


In this section we recall the definition and some important properties of $\Gamma$-convergence
and the functions of bounded variation. 
The concept of $\Gamma$-convergence was first introduced by De Giorgi in the 1970's.
For a comprehensive presentation on the topic see \cite{dalmaso}.
Let $X$ be a separable Banach space endowed with a topology $\tau$ and let 
$G,G_j:X\to [-\infty,\infty]$ for all $j\in\N$.
\begin{definition}
  \label{def:gconv}
  We say that $G_j$ $\Gamma$-converges to $G$ for the topology $\tau$
  and denote $G=\Glim_{j\to\infty} G_j$ if
  \begin{itemize}
  \item[(i)] For every $x\in X$ and for every sequence $x_j$ $\tau$-converging to $x$ in $X$ we have
      $G(x) \leq \liminf_{j\to\infty} G_j(x_j)$.
  \item[(ii)] For every $x\in X$ there exists a sequence $x_j$ $\tau$-converging to $x$ in $X$ such that
      $G(x) \geq \limsup_{j\to\infty} G_j(x_j)$.
  \end{itemize}
\end{definition}
\noindent Note that an equivalent definition is obtained by replacing condition $(ii)$ with
\begin{itemize}
\item[$(ii')$] For every $x\in X$ there exists a sequence $x_j$ $\tau$-converging to $x$ in $X$ such that
    $G(x) = \lim_{j\to\infty} G_j(x_j)$.
\end{itemize}

\begin{definition}
  A functional $G : X \to [-\infty,\infty]$ 
  is said to be {\em coercive} if condition $\lim_{j\to\infty}\norm{x_j}_X = \infty$
  implies $\lim_{j\to\infty} G(x_j) = \infty$.
  We call a sequence of functionals $G_j: X \to [-\infty,\infty]$, $j\in\N$, 
  {\em equicoercive} in topology $\tau$
  if for every $t\geq 0$ there exists
  a compact set $K_t\subset X$ such that $\{x\in X \; | \; G_j(x)\leq t\} \subset K_t$ for all $j\in\N$.
\end{definition}

The following theorem summarizes some of the known results regarding $\Gamma$-conver\-gence.
For proofs see \cite{dalmaso}.

\begin{theorem}
  \label{gconv_thm1}
  Let $G,G_j:X\to [-\infty,\infty]$, $j\in\N$, be a sequence of 
  equicoercive functionals in topology $\tau$ and $G = \Glim_{j\to\infty} G_j$.
  Then the following two properties hold:
  \begin{itemize}
  \item[(i)] If the $\Gamma$-limit of $G_j$ exists, it is unique and lower semi-continuous.
  \item[(ii)] For any continuous $H:X\to \R$ we have $G+H = \Glim_{j\to\infty} (G_j+H)$.
  \item[(iii)] Let $x_j\in X$ be such that
      $\abs{G_j(x_j)-\inf_{x\in X} G_j(x)} \leq \delta_j$
    where $\delta_j\to 0$. Then any accumulation point $y$ of $\{x_j\}_{j=1}^\infty \subset X$ is a  minimizer of $G$ and moreover $\lim_{j\to\infty}G_j(x_j) = G(y)$.
  \end{itemize}
\end{theorem}

Notice the immediate corollary to $(iii)$: suppose the assumptions 
in Theorem \ref{gconv_thm1} hold and $x_j$ is a minimizer of $G_j$ for $j\in\N$.
Then any converging subsequence of $\{x_j\}_{j=1}^\infty\subset X$ converges to a minimizer of $G$.

Let us now turn to the related function spaces.
Let $u: \T \to \R$ be a measurable function and fix $t\in\T$. We say that $z\in \R\cup \{\infty\}$ is
the {\em approximate limit} of $u$ at $t$ and write $z = \aplim_{s\to t}u(s)$
if for every neighbourhood ${\mathcal T}$ of $z$ in $\R\cup \{\infty\}$ it holds that
\[
  \lim_{\rho\to 0} \frac{1}{\rho} \abs{\{s\in\T \; | \; |s-t|<\rho, u(s)\notin {\mathcal T}\}} = 0.
\]
We use notation $\tilde u(t) = \aplim_{s\to t }u(s)$ when the limit exists.
Denote the set of points $t\in\T$ where the approximate limit does not exist by $S_u$.
When $u\in L^1(\T)$ it follows that $|S_u| = 0$, see \cite{AFP}.  

Denote by $\bv(\T)$ the Banach space of {\em functions of bounded variation}. A function $u$
belongs to $\bv(\T)$ if and only if $u\in L^1(\T)$ and its distributional derivative $Du$ is a bounded signed
measure. We endow $\bv(\T)$ with the usual norm $\norm{u}_{\bv} = \norm{u}_{L^1}+|Du|(\T)$
where $|Du|$ is the total variation of the distributional derivative.

Recall that due to the Lebesgue decomposition of measures the distributional derivative 
$Du$ can be written as a unique sum
$Du = D^a u + D^s u$ where $D^au$ is absolutely continuous and
$D^su$ is singular with respect to Lebesgue measure $|\cdot |$. Denote the density of $D^au$ with respect to the Lebesgue measure by $\nabla u$. We call function $\nabla u$ the {\em approximate gradient} of $u$.
Moreover, denote 
$D^j u = D^s u |_{S_u}$ and $D^c u = D^s u|_{\T \setminus S_u}$ where we have used notation
$\mu|_X(Y) = \mu(X\cap Y)$ for measurable sets $X,Y\subset \T$. These restrictions are
called the {\em jump part} and the {\em cantor part}, respectively.
We say $u\in \sbv(\T)$ or $u$ is a {\em special function of bounded variation} if $u \in \bv(\T)$ and $D^c u \equiv 0$.
Furthermore,  denote by $\gsbv(\T)$ the Borel functions $u:\T\to\R$ that satisfy
$\min(k,\max(u,-k)) \in \sbv(\T)$ for all $k\in\N$. The space $\gsbv$ is called the space of
{\em generalized special functions of bounded variation}.

It turns out that the generalized special functions of bounded variation inherit most important
features of $\sbv$ functions. First of all the set $S_u$ is well-defined and enumerable for $u\in \gsbv(\T)$,
and the approximate gradient $\nabla u$ exists almost every point in $\T$. We refer to \cite{AFP,braides} for
a detailed presentation on these properties.


\subsection{Mumford--Shah and Ambrosio--Tortorelli functionals}

The idea of the weak formulation of Mumford--Shah functional is to use the function
space $\gsbv$ as framework for the minimization problem and identify 
the set of jumps $K$ in \eqref{eq:strong_ms} with the set $S_u$ defined above. 
Let us drop the residual term from functional \eqref{eq:strong_ms} for the moment and denote
\[
  \MS(u,v) = \left\{
    \begin{array}{ll}
      \int_\T |\nabla u|^2dx + \sharp (S_u) & \textrm{if } u\in\gsbv(\T) \textrm{ and } v=1 \textrm{ a.e.}, \\
      \infty & \textrm{otherwise}.
    \end{array}
    \right.
\]
The role of the auxiliary function $v$ becomes clear later. The regularization term $\MS$
has been widely used in problems related to image segmentation problems. The application 
to ill-posed problems has been less extensive since in general with non-invertible forward operator $A$ 
the compactness of any minimizing sequence is not known.
For the inverse conductivity problem
in \cite{rondi, rondisantosa} the compactness is obtained by posing 
an a priori assumption that the minimizers are bounded in $L^\infty$.
In Section \ref{sec:well-posed} we prove a compactness result without such an assumption
for mildly ill-posed problems.

Next we define the Ambrosio-Tortorelli functionals \cite{AT1,AT2}.
First denote $X = H^1(\T) \times H^1(\T ; [0,1])$
and the regularizing term 
\begin{equation}
  \label{def:at}
  \at(u,v) = \int_\T \left( (\ep^2+v^2)|D u|^2 + \ep |Dv|^2 + \frac{1}{4\ep} (1-v)^2\right) dt
\end{equation}
for $(u,v) \in H^1(\T)\times H^1(\T)$. A comprehensive proof for next theorem can be found in \cite{braides} when $p=1$ and \cite{chambolle2} for the case $p=2$.

\begin{theorem}{\sc (Ambrosio-Tortorelli)}
  \label{attheorem}
  Following statement holds for $p=1$ and $p=2$.
  Define functional $\Fat :L^p(\T) \times L^p(\T) \to (-\infty,\infty]$ so that 
  \[
    \Fat(u,v) = \left\{
    \begin{array}{ll}
      \at(u,v) & \textrm{when } (u,v) \in X, \\
      \infty & \textrm{otherwise}.
    \end{array}
    \right.
  \]
  Then we have that $\Glim_{\epsilon\to 0} \Fat = \MS$ in the strong topology of $L^p(\T)\times L^p(\T)$.
\end{theorem}


\section{Main results}
\label{sec:mainresults}

Let us now return to the computational model \eqref{eq:compmodel} 
and the prior distributions introduced
in Section \ref{sec:prior_model}. For the results shown in \cite{tapio} no dependence of 
the discretization parameters $k$ and $n$ is assumed. However, in this paper we 
need to require that $k$ and $n$ are coupled, i.e., the discretization
can be characterized with only one parameter
($k=k(n)$ and $\lim_{n\to\infty} k(n) = \infty$).
For the sake of clarity in the following we assume $k=n$ and hence we drop the notation $k$.
Furthermore, the computational model \eqref{eq:compmodel} becomes simply
\begin{equation}
  \label{eq:simple_model}
  M_n = A_n U_n + \Ec_n.
\end{equation}
Before stating the assumptions concerning the Bayesian inverse problems \eqref{eq:simple_model}
for $n\in\N$ let us first introduce a definition similar to the one used in \cite{tapio,LSS}.
\begin{definition}
  \label{def:pk}
  The finite dimensional measurement projections $P_n$, $n\in\N$, are called {\em proper
  measurement projections} if they satisfy following conditions:
  \begin{itemize}
  \item[(i)] We have $\ran(P_n)\subset H^1(\T)$ and it holds that 
  $\norm{P_n}_{\Ll(H^1)}\leq C$ and $\norm{P_n}_{\Ll(L^2)}\leq C$ for some constant $C$
    with all $n\in\N$.
  \item[(ii)] For $t\in\{-1,0,1\}$ we have $\lim_{n\to\infty} \norm{P_n f-f}_{H^t} = 0$
    for all $f\in H^t(\T)$.
  \item[(iii)] For all $\phi,\psi\in L^2(\T)$ it holds that
      $\bra P_n \phi,\psi\cet_{L^2} = \bra \phi,P_n \psi\cet_{L^2}$.
  \end{itemize}
\end{definition}
For a discussion about the assumptions regarding the measurement see \cite{LSS}.
In Section \ref{sec:example} we provide an example of projections $P_n$ that satisfy Definition
\ref{def:pk}.
\begin{assumption}
  \label{ass:meas}
  For the problems in equation \eqref{eq:simple_model} 
  there exists proper measurement projections $P_n$, $n\in\N$, and fixed parameters 
  $\alpha\in \R$, $\ep>0$, and $s>0$ such that
  \begin{itemize}
  \item[(i)] there exists a bounded linear operator $A: L^2(\T)\to L^2(\T)$ 
    which satisfies
    \begin{equation}
      \label{alaraja}
      \norm{u}_{H^{-s}} \leq C \norm{Au}_{L^2}
    \end{equation}
    for any $u\in L^2(\T)$ with some constant $C>0$ 
    and
    $A_n = P_n A$ for all $n\in\N$.
  \item[(ii)] The additive noise $\Ec_n$ is a Gaussian random variable in $PL(n)$ such that
    $\expec \Ec_n = 0$ and for any $\phi,\psi \in L^2(\T)$ the covariance satisfies 
    \[
      \expec \bra \Ec_n,\phi\cet_{L^2} \bra \Ec_n, \psi\cet_{L^2}
      = N^{-\alpha} \bra P_n \phi,P_n \psi\cet_{L^2}.
    \]
  \item[(iii)] The prior structure is modelled 
    with random variables $(U_{n,\ep}^\alpha, V_{n,\ep}^\alpha)$ 
    introduced in Section \ref{sec:prior_model}.
  \item[(iv)] The measurements $m,m_n\in L^2(\T)$, $n\in \N$, satisfy $\lim_{n\to\infty} m_n = m$ in $L^2(\T)$.
  \end{itemize}
\end{assumption}
Notice that condition $(ii)$ means that $\Ec_n$ has white noise statistics and
in case $\alpha=0$ 
the random variables $\Ec_n$ 
convergence to white noise in the sense of generalized random variables as $n\to\infty$
\cite{LSS,rozanov}.
Now with Assumption \ref{ass:meas} the variational problem of finding the 
MAP estimates for equation \eqref{eq:simple_model} becomes
\begin{equation}
  \label{def:map_est}
  {\rm min}_{(u,v)\in PL(n)\times PL(n)} 
  \left(\Fda(u,v) + \norm{P_n A u - m_n}^2_{L^2}\right).
\end{equation}

Below we study the behaviour of the MAP estimates with respect to parameters
$n,\alpha$ and $\epsilon$ using the 
variational approximation methods presented in Section
\ref{sec:Glim_and_BV}. In order to describe the $\Gamma$-limits of the functionals in equation
\eqref{def:map_est} we have to introduce some new notations.

Let us first denote an auxiliary domain
\begin{equation}
  \label{Xep_domain}
  X_\ep = H^1(\T) \times H^1(\T; [0,1+30\ep])
\end{equation}
for sufficiently small $\ep$.
Details about this choice of domain are given in Appendix \ref{sec:choiceofa}
and we discuss it in more detail below.
For now, it suffices to point out that 
the domain $X_\ep$ formally approaches
$X$ when $\ep$ decreases.
Denote the auxiliary operators
\[
  \pone(v) = \int_\T - \log(\ep^2+v^2) dt 
\]
and
\[
  \ptwo(u,v) = \int_\T (\ep^2+v^2)(2 \ep^q(Qu) Du+ \ep^{2q}(Qu)^2) dt 
\]
for $(u,v)\in H^1(\T)\times H^1(\T)$. We motivate these notations after the next definition. 
Recall that $X$ denotes the domain $H^1(\T)\times H^1(\T; [0,1])$.
\begin{definition}
  \label{def:fe_and_fet}
  Let us define functionals $\Fbe : L^1(\T)\times L^1(\T) \to (-\infty,\infty]$, $\ep>0$ and 
  $\alpha\geq 1$, so that for $\alpha=1$
  \[
    \Fe^1(u,v) = \left\{
    \begin{array}{ll}
      \pone(u,v)+\at(u,v)+\ptwo(u,v) & \textrm{when } (u,v) \in X_\ep, \\
      \infty & \textrm{otherwise}.
    \end{array}
    \right.
  \]
  and for $\alpha>1$
  \[
    \Fbe(u,v) =  \left\{
    \begin{array}{ll}
      \at(u,v)+\ptwo(u,v) & \textrm{when } (u,v) \in X, \\
      \infty & \textrm{otherwise}.
    \end{array}
    \right. 
  \]
\end{definition}

Let us now discuss Definition \ref{def:fe_and_fet}.
The reason for the particular choice of $X_\ep$ is two-fold. First, it turns out that
the minimizer of functional $\pone+\at+\ptwo$ in $H^1(\T)\times H^1(\T)$ may be located outside $X$.
Secondly, a pointwise bound for $v$ provides easier proofs
concerning the $\Gamma$-convergence results of
functionals $\Fbe$ in Section \ref{sec:convergence}. 

Furthermore, it is straightforward to see that
\begin{equation}
 \label{aux_form_fe}
  \at(u,v)+\ptwo(u,v) = \int_\T \left((\ep^2+v^2)|\td u|^2 + \ep |Dv|^2 + \frac{1}{4\ep} (1-v)^2\right) dx
\end{equation}
everywhere in $H^1(\T)\times H^1(\T)$. 
Hence the role of $\ptwo$ can be understood as a small perturbation that 
yields a lower bound for $|Qu|$ and thus coersivity for $\Fbe$. 
On the other hand, compared to the Ambrosio-Tortorelli approach, a new term
$\pone$ appears due to the Bayesian hierarchical modelling. 

In addition to problem \eqref{def:map_est}, we will consider three 
different minimization problems throughout the paper. 
Two of them are the
modified Ambrosio--Tortorelli minimization problem
\begin{equation}
  \label{eq:fbe+res}
  \min_{(u,v)\in H^1(\T)\times H^1(\T)} \Fbe(u,v) + \norm{Au-m}^2_{L^2}
\end{equation}
and the Mumford--Shah problem
\begin{equation}
  \label{eq:ms+res}
  \min_{(u,v)\in L^1(\T)\times L^1(\T)} \MS(u,v) + R(u).
\end{equation}
In \eqref{eq:ms+res} we assume that $A: L^p(\T) \to L^p(\T)$ is continuous for $p\in \{1,2\}$
and the residual $R(u)$ is defined by
\[
  R(u) = \left\{
    \begin{array}{ll}
      \norm{Au-m}^2_{L^2} & \textrm{when } Au\in L^2(\T), \\
      \infty & \textrm{otherwise}.
    \end{array}
    \right. 
\]
In the following we often use notation $\norm{Au-m}^2_{L^2}$ for $R(u)$
when convenient.
To describe cases when the edge-preserving property of MAP estimates is lost asymptotically we consider
the Tikhonov-type minimization problem
\begin{equation}
  \label{tik}
  H(u) = \left\{
    \begin{array}{ll}
      \int_\T |Du|^2dt + \norm{Au-m}^2_{L^2} & \textrm{when } u\in H^1(\T), \\
      \infty & \textrm{otherwise}.
    \end{array}
    \right. 
\end{equation}
Notice that the solution to problem \eqref{tik} is obtained by $u_{min} = (-\Delta  + A^* A)^{-1}A^*m$.
With the definitions given above we are ready to state the main results.
Denote the conditional mean estimates introduced in Section \ref{sec:bayes} of problem
\eqref{eq:simple_model} by $(u^{CM}_{n,\ep},v^{CM}_{n,\ep})$.

\begin{theorem}
  \label{neg_result}
  Let the computational model \eqref{eq:simple_model} satisfy Assumption \ref{ass:meas}
  with prior parameters  $s>0$, $\alpha = 0$, and $\ep>0$ and let
  the operator $A:L^2(\T)\to H^1(\T)$ be bounded. Then the following statements hold:
  \begin{itemize}
    \item[(i)] The CM estimates $(u^{CM}_{n,\ep},v^{CM}_{n,\ep})$ 
    	converge in $L^2(\T)\times L^2(\T)$ as $n\to \infty$.
    \item[(ii)] The MAP estimates $(u^{MAP}_{n,\ep},v^{MAP}_{n,\ep})$ diverge as $n\to \infty$.
    \end{itemize}
    In addition, the following holds for coupled parameters:
    \begin{itemize}
    \item[(iii)]
  If $\ep=\ep(n)\to 0$ as $n\to\infty$ then
  it follows that either the minimum values in formula \eqref{def:map_est} diverge to $-\infty$
  or the MAP estimates $(u^{MAP}_{n,\ep(n)},v^{MAP}_{n,\ep(n)})$ converge towards a minimizer of 
  the functional \eqref{tik}.
  \end{itemize}
\end{theorem}

The statement $(i)$ of Theorem \ref{neg_result} is proved in \cite{tapio} and 
statements $(ii)$ and $(iii)$ are proven in Section
\ref{sec:non-edge-preserving}. Notice that even the coupling of $\ep$ and $n$ in statement $(iii)$
do not yield convergence to Mumford--Shah minimizers. 
Namely, the diverging minimum values immediately contradict with condition $(i)$ in Definition \ref{def:gconv}
since functional $MS$ is positive. Furthermore, the convergence to a minimizer of functional \eqref{tik} implies
that the edge-preserving property of the MAP estimates is lost. We point out that statement $(iii)$ does not imply that
this property is lost with {\em all couplings}. 

Our main positive result regarding the convergence of the MAP estimates is the following.
\begin{theorem}
  \label{pos_result}
  Let the computational model \eqref{eq:simple_model} satisfy Assumption \ref{ass:meas}
  with prior parameters $s<\frac 12$, $\alpha \geq 1$, and $\ep>0$
  and let the operator $A:L^p(\T)\to L^p(\T)$ be bounded for $p=1,2$.
  Then 
    \begin{itemize}
  \item[(i)] The MAP estimates $(u^{MAP}_{n,\ep},v^{MAP}_{n,\ep})$ 
  have a subsequence converging to a minimizer
    $(u_{\ep},v_{\ep})$ 
    of problem \eqref{eq:fbe+res} in the weak topology of $H^1(\T)\times H^1(\T)$ as $n\to\infty$.
  \item[(ii)] The minimizers $(u_{\ep},v_{\ep})\in H^1(\T)\times H^1(\T)$, $\ep>0$, 
  of the problem \eqref{eq:fbe+res}
  have a subsequence converging
  to a minimizer of the Mumford--Shah problem \eqref{eq:ms+res} in 
  $L^1(\T)\times L^1(\T)$ as $\ep\to 0$. 
  \end{itemize}
\end{theorem}

The result $(ii)$ in Theorem \ref{pos_result} can be also considered as a new interpretation
of the Mumford--Shah functional; the minimizer of the Mumford--Shah functional can be approximated
by the MAP estimates of Bayesian inverse problems.
The proof for Theorem \ref{pos_result} is given in Section \ref{sec:convergence}.


\section{Well-posedness of the minimization problems}
\label{sec:well-posed}

In this section we study the properties of the individual problems \eqref{def:map_est}, \eqref{eq:fbe+res}
and \eqref{eq:ms+res} with fixed parameters $\ep,\alpha$ and $n$. 
Our aim is to show three results.
First, the existence of a minimizer of problem \eqref{eq:ms+res} is proved in
Theorem \ref{MSexistence}.
Second, we show in Lemma \ref{inf_domain} that with the choice of domain $X_\ep$ 
we do not exclude any pairs $(u,v)\in H^1(\T)\times H^1(\T)$ which give a smaller value in the problem \eqref{eq:fbe+res}.
Finally, we show that functionals $\Fda$ and $\Fbe$ are coercive in $H^1(\T)\times H^1(\T)$
which yields the existence of minimizers in problems \eqref{def:map_est} and \eqref{eq:fbe+res}.

Let us now study the existence of solution to problem \eqref{eq:ms+res}. 
The following compactness and semi-continuity theorem in $\gsbv$ is well-known. 
\begin{lemma}
  \label{ms_existence}
  Suppose a sequence $\{u_j\}_{j=1}^\infty \subset \gsbv(\T)$ satisfies
  \begin{equation}
    \label{cpt_condition}
    \norm{u_j}_{L^p}+\sharp (S_{u_j}) + \int_\T |\nabla u_j|^2dt \leq C
  \end{equation}
  for some $1<p\leq 2$.
  Then there exists $u\in \gsbv(\T)\cap L^p(\T)$ and a subsequence $\{u_{j_k}\}_{k=1}^\infty$ such that
  \begin{itemize}
  \item[(i)] $u_{j_k} \to u$ strongly in $L^1(\T)$,
  \item[(ii)] $\nabla u_{j_k} \weakly \nabla u$ weakly in $L^2(\T)$ and
  \item[(iii)] $\sharp (S_u) \leq \liminf_{k\to\infty} \sharp (S_{u_{j_k}})$.
  \end{itemize}
\end{lemma}

\begin{proof}
  If $\{u_j\}_{j=1}^\infty\subset \gsbv(\T)$ satisfies \eqref{cpt_condition}
  then also
  \begin{equation}
    \label{aux_cpt_condition}
    \norm{u_j}_{L^p}+\sharp (S_{u_j}) + \int_\T |\nabla u_j|^pdt \leq C.
  \end{equation}
  By \cite{a1,a2,a3} (see \cite[Thm. 2.1.]{cortesani} for a short exposition)
  it holds that there exists a subsequence $\{u_{j_k}\}_{k=1}^\infty$ and $u \in \gsbv(\T)\cap L^p(\T)$
  such that conditions $(i)$ and $(iii)$ are satisfied. Furthermore, since $\nabla u_j$
  is bounded in $L^2(\T)$ due to the Banach--Alaoglu
  theorem we can extract a subsequence such that condition $(ii)$ holds.
\end{proof}

Next lemma shows that the assumption in Lemma \ref{ms_existence} is
in a sense self-improving and one can extend it
for the purpose of mildly ill-posed problems.

\begin{lemma}
  If a function $u \in \gsbv(\T)\cap L^1(\T)$ satisfies
  \begin{equation}
    \label{secondary_cpt_condition}
    \norm{u}_{H^{-s}}+\sharp (S_{u}) + \int_\T |\nabla u|^2dt \leq C'
  \end{equation}
  for some $0\leq s<\frac 12$ then for $p>1$ such that $\frac 1p = s + \frac 12$
  it also satisfies inequality \eqref{cpt_condition} for some constant $C$ depending
  on $s$ and $C'$.
\end{lemma}

\begin{proof}
  Let us denote by $t_j$ the points in $S_{u}$,
  such that $S_{u} = \{t_1, t_2, ..., t_{L}\}$, where $t_1<t_2<...<t_L$ and $L = \sharp (S_{u})$ 
  is bounded.
  Furthermore, denote by $I_j = (t_{j-1},t_j)\subset \T$, $1\leq j\leq L$, 
  the interval between neighboring points.
  Here $t_0$ and $t_{L}$ were identified. 
  We can estimate the average of $u$ over interval $I_j$ by
  \[
    \frac{1}{|I_j|}\abs{\int_{I_j} u_k dt} 
    \leq C|I_j|^{-\frac 12-s} \norm{u}_{H^{-s}(I_j)},
  \]
  where we have used Lemma \ref{dual_bound}. Now the Poincare inequality 
  states that
  \[
    \norm{u-\frac{1}{|I_j|}\int_{I_j}u dt}_{L^p(I_j)} \leq C |I_j| \norm{\nabla u}_{L^p(I_j)}
  \]
  and we obtain $\norm{u}_{L^p(I_j)} \leq C(|I_j|+1)$.
   By using the knowledge 
   $\sum_{j=1}^{L} |I_j| = 1$ we deduce that $\norm{u}_{L^p(\T)} \leq C''$
   where the constant $C''$ depends only on $s$ and $C'$.
   This proves the claim.
\end{proof}

Clearly any sequence $u_j$ satisfying inequality \eqref{cpt_condition}
belongs to $L^\infty(\T)$ and thus also $\sbv(\T)$.
However the bound in \eqref{cpt_condition} does not control this norm and hence
without any additional bound in $L^\infty$ the limit does not necessarily belong to $\sbv(\T)$.
As the existence of a Mumford--Shah minimizer has interest
for inverse problems in general we have formulated an independent proof
to the following theorem.

\begin{theorem}
  \label{MSexistence}
  Let $A$ be a bounded linear operator in $L^p(\T)$ for $p=1,2$ such that 
  it satisfies inequality \eqref{alaraja} for some $s<\frac 12$.
  Then the minimization problem
  \begin{equation}
    \label{infms}
    \inf_{u\in L^1(\T)}\left(\MS(u,\onef) + R(u)\right)
  \end{equation}
  has a solution $u\in \gsbv(\T)\cap L^1(\T)$.
\end{theorem}

\begin{proof}
  Assume that $\{u_j\}_{j=1}^\infty \subset L^1(\T)$
  is a minimizing sequence, i.e.,
  \[
     \inf_{u\in L^1(\T)}\left(\MS(u,\onef) + R(u)\right)  = 
    \lim_{j\to\infty} \left(\MS(u_j,\onef) + R(u_j)\right).
  \]
  Then the sequence $u_j$ satisfies inequality \eqref{secondary_cpt_condition} which in turn
  yields conditions in Lemma \ref{ms_existence}. 
  Consequently, we may extract a 
  subsequence $u_{j_k}$ converging in $L^1(\T)$ to
  $u\in \gsbv(\T)\cap L^1(\T)$. Notice that the residual term $R(u)$ is lower semicontinuous
  in the $L^1(\T)$ topology.
  Denoting the infimum in \eqref{infms} by ${\mathcal I}$ we obtain using Lemma \ref{ms_existence} that
  \begin{eqnarray*}
     {\mathcal I} & \leq & \MS(u,\onef) + R(u) \\
     & \leq & \liminf_{k\to\infty} \MS(u_{j_k},\onef) + \liminf_{k\to\infty} R(u_{j_k}) \\
    & \leq & \liminf_{k\to\infty} \left(\MS(u_{j_k},\onef) + R(u_{j_k})\right) \leq {\mathcal I}.
  \end{eqnarray*}
  The claim follows from $(u,\onef)$ being a minimizer.
\end{proof}

Next we discuss the choice of domains $X$ and $X_\ep$ in Definition \ref{def:fe_and_fet}.
Denote by $\psi_r, \Psi_r : \R \to [0,r]$, $r>0$, 
the functions
$\psi_r(t) = (r-|r-t|) \chi_{[0,2r]}(t)$
and
\begin{equation}
  \label{psifunction}
  \Psi_r(t) = \sum_{j\in \Z} \psi_r(t-2jr).
\end{equation}
We notice that for any function $f$ and $r>0$ the mapping $\Psi_r\circ f$ satisfies
$0\leq (\Psi_r\circ f)(t) \leq r$ for all $t\in\T$. Due to such property we call this operation {\em folding}.
We list the following three properties of $\Psi_r$ as a lemma.
\begin{lemma}
  \label{psilemma}
  For any $f\in H^1(\T)$ it holds that
  \begin{itemize}
    \item[(i)] $|(\Psi_r \circ f)(t)| \leq |f(t)|$ for any $t\in\T$,
    \item[(ii)] $|r-(\Psi_r \circ f)(t)| \leq |r-f(t)|$ and
    \item[(iii)] $|D (\Psi_r \circ f)| = |D f|$ almost everywhere on $\T$.
  \end{itemize}
\end{lemma}
\begin{proof}
  The first claim is obvious since $\Psi_r(f(t)) = {\rm sgn}(f(t))|f(t)|$
  when $f(t) \in [-r,r]$ and also $0\leq \Psi_r(f(t))\leq r$ for any $t\in\T$.
  Claim $(ii)$ also follows from the definition of function $\Psi_r$.
  For claim $(iii)$, notice that since $f\in H^1(\T)$ by the Sobolev embedding theorem $f$ must be bounded, i.e.,
  $\sup_{t\in\T} |f(t)| < C$. In consequence, $\Psi_r\circ f$ can be written as a finite sum over
  functions $\psi_r(f(\cdot)-2jr)$. Now the result follows from a generalization of the chain rule
  (see e.g. \cite{adm}).
\end{proof}

In particular, Lemma \ref{psilemma} yields that $\Psi_r\circ f \in H^1(\T)$ whenever $f\in H^1(\T)$.
In the proof of next lemma we use the idea that in some cases $\Psi_r\circ v$
with suitable choices of $r>0$ produces lower value than $v$ for the considered functionals.
Consequently, we obtain information that the minimizers must lie in $X$ or $X_\ep$.
\begin{lemma}
  \label{inf_domain}
  Let $\alpha\geq 1$.
  For every $v\in H^1(\T)$ there exists $w\in H^1(\T ; [0,1+30\ep])$ such that
  \begin{equation}
    \label{fbe_domain}
    \Fbe(u,w) \leq 
    \delta_{\alpha,1}\pone(v)+\at(u,v) + \ptwo(u,v).
  \end{equation}
  for all $u\in H^1(\T)$ where $\delta_{\alpha,1}=1$ when $\alpha=1$ and is otherwise zero.
\end{lemma}

\begin{proof}
  Consider function $\Psi_1$ defined by equation \eqref{psifunction}. 
  Due to Lemma \ref{psilemma} and equation \eqref{aux_form_fe} we have
  $\Fe^\alpha(u,\Psi_1\circ v) \leq \Fe^\alpha(u,v)$
  for any $(u,v)\in H^1(\T)\times H^1(\T)$ and $\alpha>1$.  
  This immediately yields inequality \eqref{fbe_domain} with $w=\Psi_1\circ v$ for $\alpha>1$.
  
  Let us then consider the case $\alpha=1$ and let $(u,v) \in H^1(\T)\times H^1(\T)$. 
  To apply folding denote $E_-=\{v(t)<0\}$,
  $E_0 = \{0\leq v(t) \leq 1+30\ep\}$ and $E_+ = \{v(t)>1+30\ep\}$ and by
  $\onef_E$ the indicator function of $E$. We write
   \[
    v = v \cdot \onef_{E_-} + v \cdot \onef_{E_0} + v \cdot \onef_{E_+} = v_- + v_0 + v_+.
  \]
  We construct $w$ by applying the folding operation to each restriction separately.
  First, recall identity \eqref{aux_form_fe} and denote
  \begin{equation}
    \label{integrand}
    G_\ep(v,E) = \int_E \left(-\log(\ep^2+v^2)+\frac{1}{4\ep}(1-v)^2\right)dt.
  \end{equation}
  with any measurable $E\subset \T$.
  Denote then 
  \[
    w_-=\Psi_1\circ v_-, \quad w_+=\Psi_{30\ep}\circ (v_+-1) +1  \quad
    {\rm and} \quad w = w_- + v_0 + w_+.
  \]
  Clearly $0\leq w \leq 1+30\ep$ and $(u,w)\in X_\ep$.

  First, we see due to Lemma \ref{psilemma} claim $(i)$ that 
  \[
  |w(t)| = |w_-(t)|+|v_0(t)|+|w_+(t)|\leq |v(t)|
  \]
  for all $t\in\T$.
  Furthermore, claim $(iii)$ in Lemma \ref{psilemma} implies $|Dw(t)|=|Dv(t)|$ almost everywhere on $\T$. 
  These yield
  \begin{equation}
    \label{fbe_dom_aux1}
    \int_\T ((\ep^2+w^2)|\td u|^2 + \ep|Dw|^2)dt \leq
    \int_\T ((\ep^2+v^2)|\td u|^2 + \ep|Dv|^2)dt.
  \end{equation}
  
  Let us next consider the integrand
  $g_\ep(t) = -\log(\ep^2+t^2)+\frac{1}{4\ep}(1-t)^2$
  in equation \eqref{integrand} and apply the
  results in Lemma \ref{aux_function} to functions $w_-$, $w_0$ and $w_+$.
  Due to claim $(ii)$ and $(iii)$ in Lemma \ref{aux_function}
  it is straightforward to see
  $G(w_-,E_-)\leq G(v_-,E_-)$.
  Furthermore, the claim $(i)$ in Lemma \ref{aux_function}
  implies
  $G(w_+,E_+)\leq G(v_+,E_+)$.
  From this we conclude
  \begin{eqnarray}
    \label{fbe_dom_aux2}
    G(w,\T) & = & G(w_+,E_+)+G(v_0,E_0)+G(w_-,E_-)\nonumber \\
    & \leq & G(v_+,E_+) + G(v_0,E_0)+G(v_-,E_-) = G(v,\T).
  \end{eqnarray}
  Now inequalities \eqref{fbe_dom_aux1} and \eqref{fbe_dom_aux2} together with
  identity \eqref{aux_form_fe} yield the result.
\end{proof}

\begin{theorem}
  The functionals $\Fbe$ for $\alpha\geq 1$ and $\Fda$ for $\alpha\in\R$ are coercive in $H^1(\T)\times H^1(\T)$ for any 
  fixed $\alpha$, $n\in\N$ and $\ep>0$.
\end{theorem}

\begin{proof}
  Recall that a functional $G: X\to \R$ is coercive if we have a lower bound $G(x)\geq C\norm{x}_X$ for $x\in X$ such
  that $\norm{x}_X$ is large enough.
  By the Lemma \ref{infFep} in Appendix we know that the functionals are bounded from below.  
  One can deduce that 
  \[
    \int_\T \ep^2 |\td u|^2dt= \int_\T \left(\ep^2 |Du|^2dt + \ep^{2+2q}(Qu)^2\right)dt \geq C(\ep) \norm{u}^2_{H^1}.
  \]
  The lower bound for $\norm{v}^2_{H^1}$ can be obtained from the term 
  $\int_\T \left(\ep |Dv|^2+\frac{1}{4\ep}(1-v)^2\right)dt$. Hence it follows that both
  $\Fbe(u,v)$ and $\Fda(u,v)$ go to infinity when $\norm{u}_{H^1}$ or $\norm{v}_{H^1}$ goes to infinity.
\end{proof}


\section{Non-edge-preserving scaling}
\label{sec:non-edge-preserving}

In this section we study the case when $s>0$ and  $\alpha=0$ and prove Theorem \ref{neg_result}. Recall 
that the claim $(i)$ is shown in \cite{tapio}.

\vspace{0.5cm}
{\em Proof of claim $(ii)$ in Theorem \ref{neg_result}.}
Consider the value of $F^0_{\ep,n}$ at function $(u(t),v(t))\equiv (0,s)$
where $s>1$, namely,
\[
  F^0_{\ep,n}(u,v) = g_N(s) = -N\log(\ep^2+s^2) + \frac{1}{4\ep}(1-s)^2
\]
where $N=2^n$.
With fixed $s>1$ we have then $\lim_{n\to\infty} g_N(s) = -\infty$.
Also, it is easy to see that the minimizing values $s = s(N)$ of $g_N(s)$
go to infinity if $n\to\infty$. 
Suppose now that the pair $(u_n,v_n)\in PL(n)\times PL(n)$ is a minimizer for problem
\eqref{def:map_est} and that $(u_n,v_n)$ converges in $H^1(\T)\times H^1(\T)$.
Clearly then
\[
  F^0_{\ep,n} (u_n,v_n)+\norm{A_nu_n-m_n}_{L^2}^2 \leq g_N(s(N)) + \norm{m_n}^2_{L^2}
\]
for all $n\in\N$. Since the terms of $F^0_{\ep,n}$ are all positive except for the logarithm
term and the measurements $m_n$ are bounded in $L^2(\T)$ we have
\begin{equation}
  \label{aux_ineq_nonedge1}
  \int_\T -N \log(\ep^2+(Q_nv_n)^2)dt \leq g_N(s(N)) + C
\end{equation}
for some constant $C>0$. The assumption that $v_n$ converge in $H^1(\T)$ yields that $\norm{v_n}_\infty<C'$
for all $n\in\N$ with some $C'>0$. Thus, inequality \eqref{aux_ineq_nonedge1} implies
$ - N C'' \leq g_N(s(N)) + C$
for all $n\in\N$ where $C''=\log(\ep^2+(C')^2)$. This leads to a contradiction since
$\lim_{n\to\infty} \frac{N}{g_N(s(N))} = 0$
and proves the claim.
\hfill $\Box$

\vspace{0.5cm}
The immediate question after the result $(ii)$ in Theorem \ref{neg_result} is
whether an appropriate coupling of $\ep$ and $n$ guarantees the convergence of MAP estimates.
In the following we give some negative results about this.

Consider first how the discretization scheme affects the convergence in $\norm{\cdot}_\infty$-norm.

\begin{theorem}
  \label{vn_uniform}
  Assume $u_n,v_n\in PL(n)$ for $n\in\N$, $N=2^n$ and $F_{\ep,n}^0(u_n,v_n) \leq C$
  for some constant $C>0$. Then there exists a constant $C'>0$ such that
  \[
    \norm{1-v_n}_\infty \leq C' \sqrt{\ep N + \ep^2 N^3}.
  \]
\end{theorem}

\begin{proof}
  The boundedness of $F_{\ep,n}^0(u_n,v_n)$ and Lemma \ref{infFep} yield
  \[
    -C_1\ep N^2 + \int_\T \frac{1}{\ep} (1-v_n)^2 dt \leq C_2
  \]
  for some constants $C_1,C_2>0$. This immediately results to
  \begin{equation}
    \label{vn_unif_aux}
    \norm{1-v_n}^2_{L^2} \leq C_2\ep+ C_1\ep^2N^2.
  \end{equation}
  First denote $t_k = k/N$ for all $0\leq k<N$.
  Suppose that $f\in PL(n)$ achieves its maximum at $t_j$ and denote by $\phi_j\in PL(n)$ a function
  that satisfies $\phi_j(t_k) = \delta_{jk}$ for all $0\leq k< N$. Then by using the simple fact
  that $\int_{t_{j-1}}^{t_{j+1}} f(t_j)\phi_j(t) dt \leq \int_{t_{j-1}}^{t_{j+1}} |f(t)|dt$ we have
  \[
    \norm{1-v_n}_\infty = N \int_{t_{j-1}}^{t_{j+1}} \norm{1-v_n}_\infty \phi_j(t)dt 
    \leq \sqrt{N} \sqrt{\int_\T |1-v_n(t)|^2dt}    
  \]
  where in the last inequality we have used the Cauchy-Schwarz inequality.
  Now inequality \eqref{vn_unif_aux} proves the claim.
\end{proof}

\begin{corollary}
  \label{nogconvergence}
  Let $\ep=\ep(n)$ such that $\lim_{n\to \infty} \ep(n) = 0$ and
  let $(u_n,v_n)\in PL(n)\times PL(n)$ be a minimizer of $F^\alpha_{\ep,n}$. Then
  the following statements hold:
  \begin{itemize}
    \item[(i)] If $\lim_{n\to\infty} \sqrt{\ep(n)} 2^n < \infty$ then 
      the function $v_n$ converges uniformly to $\onef$ with respect to $n$.
    \item[(ii)] If $\lim_{n\to\infty} \sqrt{\ep(n)} 2^n = \infty$ 
      then the minimum values $F^0_{\ep,n}$ diverge, i.e.,
     \[
     \lim_{n\to\infty} F^0_{\ep,n} (u_n,v_n) = -\infty.
     \]    
  \end{itemize}
\end{corollary}

\begin{proof}
Let us first notice that 
\[  
  F^0_{\ep,n}(u_n,v_n) \leq \inf_{v \in PL(n)} F^0_{\ep,n}(0,v) = \inf_{v\in PL(n)}
  \int_\T \left(-2^n \log(\ep^2+(Q_nv)^2) + \frac{1}{4\ep}(1-v)^2\right)dt.
\]
The statement $(ii)$ follows the upper bound given in Lemma \ref{infFep}.

Assume now that $\lim_{n\to\infty} \sqrt{\ep(n)} 2^n < \infty$. 
By using inequality $\log(1+x) \leq x$ it follows that also
$\lim_{\ep\to 0} \left(2^n \log(1+O(\sqrt{\ep}))\right) < \infty$.
By simple computations one can show $\lim_{\ep\to 0} \frac{\ep^p}{\log(1+O(\sqrt{\ep}))} = 0$
for any $p>\frac 12$ and hence the quantity
\[
  2^{3n} \ep^2 = \left(2^n\log(1+O(\sqrt{\ep})) \right)^3
  \left(\frac{\ep^{\frac 2 3}}{\log(1+O(\sqrt{\ep}))}\right)^3
\]
converges to zero. The convergence of $2^{2n}\ep^2$ to zero follows by the same argument.
Consequently, the result $(i)$ follows from Theorem \ref{vn_uniform}.
\end{proof}

For $n\in\N$ define functionals
\begin{equation}
  \label{def:hn}
  H_n(u) = \int_\T h_n^2 |\td u|^2 dt + \norm{P_nA u- m_n}^2_{L^2}
\end{equation}
for $u\in PL(n)$ where $h_n\in H^1(\T; \R_+)$ converges to $\onef(t) \equiv 1$ uniformly.

\begin{theorem}
  \label{glim_hn}
  Let $\ep=\ep(n)$ such that $\lim_{n\to\infty} \ep(n) = 0$. We have that
  \begin{itemize}
  \item[(a)] $H = \Glim_{n\to\infty} H_n$ in the weak topology of $H^1(\T)$
  \item[(b)] the functionals $\{H_n\}_{n\in\N}$ are equicoercive in the 
    weak topology of $H^1(\T)$.
  \end{itemize}
\end{theorem}

\begin{proof}
  Let $u_n\weakly u$ weakly in $H^1(\T)$ as $n\to\infty$ where $u_n\in PL(n)$. 
  By lower semicontinuity of the norm we have
  \[
    \int_\T |Du|^2dt \leq \liminf_{n\to\infty} \int_\T |Du_n|^2dt  
    	= \liminf_{n\to\infty} \int_\T h_n^2 |\td u_n|^2dt.
  \]
  Furthermore, by the Sobolev embedding theorem we see that $u_n\to u$ in $L^2(\T)$ and hence
  $\lim_{n\to\infty} P_nAu_n - m_n = Au-m$
  in $L^2(\T)$. Together these imply $H(u) \leq \liminf_{n\to\infty} H_n(u_n)$.
  This proves $(i)$ in Definition \ref{def:gconv}.
  
  To prove the condition $(ii)$ in Definition \ref{def:gconv} it is sufficient to consider any sequence $u_n\in PL(n)$
  such that $u_n\to u$ in the $H^1(\T)$-norm as $n\to\infty$.
  This proves the claim $(a)$ here. Let us then study claim $(b)$ and assume that $u_n\in PL(n)$ for every $n\in\N$
  and 
  \begin{equation}
    \label{aux_hn_bound}
    H_n(u_n) = \int_\T h_n^2 |\td u_n|^2 dt + \norm{P_nA u_n- m_n}^2_{L^2} \leq C
  \end{equation}
  for some constant $C>0$. In particular, we have $\norm{P_nAu_n}_{L^2} \leq C$ for all $n\in\N$.
  
  Next we show that also the sequence $\norm{Au_n}_{L^2}$ is uniformly bounded.
  Assume for the moment that this is not the case and 
  $\lim_{n\to\infty} \norm{Au_n}_{L^2}=\infty$. Recall that operator $Q$ was defined as $Qf = (\int_\T f(t)dt) \onef$
  for any $f\in L^2(\T)$.
  Due to the inequality \eqref{aux_hn_bound} we have that $\norm{Du_n}_{L^2}\leq C$
  and, in consequence, $\norm{(I-Q)u_n}_{L^2} \leq C$.
  Moreover, this yields
  \begin{equation}
  	\label{helpeq1}
    \lim_{n\to\infty} \norm{AQu_n}_{L^2} \geq \lim_{n\to\infty} (\norm{Au_n}_{L^2} - \norm{A(I-Q)u_n}_{L^2}) = \infty.
  \end{equation}
  By setting $c_n = \int_\T u_n(t)dt$ equation \eqref{helpeq1} implies that
  \begin{equation}
    \label{eq:hnaux1}
    \lim_{n\to\infty} |c_n| \norm{A\onef}_{L^2} = \infty.
  \end{equation}
  The boundedness of $\norm{A(I-Q)u_n}_{L^2}$ yields together with $\norm{P_nAu_n}_{L^2}\leq C$ that 
  \begin{equation}
    \label{eq:hnaux2}
    |c_n| \norm{P_nA \onef}_{L^2} = \norm{P_nA Q u_n}_{L^2}\leq \norm{P_nA(I-Q)u_n}_{L^2} +\norm{P_nAu_n}_{L^2} \leq C.
  \end{equation}
  By the condition \eqref{eq:hnaux1} we have $\lim_{n\to\infty} |c_n| = \infty$ and by condition \eqref{eq:hnaux2}
  it follows that $\lim_{n\to\infty} \norm{P_nA \onef}_{L^2} = 0$.
  Due to the condition $(ii)$ 
  in Definition \ref{def:pk} 
  this yields 
  $A\onef = 0$. However, this contradicts with equation \eqref{eq:hnaux1}.
  Consequently, we have proven that $\norm{Au_n}_{L^2} \leq C$ for some constant $C>0$.

  By the assumption on $A$ with $s>0$,  we have $\|u_n\|_{H^{-s}} \leq C\norm{A u_n}_{L^2}$.
 As $1\in (H^{-s}(\T))'=H^{s}(\T)$, we have $|Qu_n| \leq C\norm{u_n}_{H^{-s}}$ and hence we obtain by the Poincare inequality 
  that $\norm{u_n}_{H^1}$ is bounded. 
  By the Banach-Alaoglu theorem there exists a converging subsequence which completes the proof.
\end{proof}

Finally we conclude this section by completing the proof of Theorem \ref{neg_result}.
\vspace{0.2cm}

{\em Proof to claim (iii) in Theorem \ref{neg_result}}. 
Suppose that $\lim_{n\to\infty} \ep(n) = 0$ and that
$(u_n,v_n)\in PL(n)\times PL(n)$ minimizes $F^0_{\ep(n),n}$.
By Corollary \ref{nogconvergence} either the minimum values of $F^0_{\ep(n),n}$ diverge to $-\infty$
or $v_n\to \onef$ uniformly. 

Consider the latter case. Then the functions $u_n$ clearly solve minimization
problems $\min_{u\in PL(n)} H_n(u)$ where $H_n(u)$ is defined in equation \eqref{def:hn}
with $h_n = \ep^2+v_n^2$. By Theorems \ref{gconv_thm1} and \ref{glim_hn} it follows that $u_n$
converge to a minimizer of $H$.
\hfill $\Box$


\section{Convergence proofs}
\label{sec:convergence}


\subsection{Convergence with respect to $n$}

In this section we prove $\Gamma$-convergence of $\Fda$ with respect to $n$ for all scalings $\alpha\geq 1$.  Throughout the section, $0<s<\frac 12$.

\begin{theorem}
  \label{gconv_for_n}
  For $\alpha\geq 1$ we have $\Fbe = \Glim_{n\to\infty} \Fda$ in the weak topology of 
  $H^1(\T)\times H^1(\T)$.
\end{theorem}

\begin{proof}
  Let us assume that $(u_n,v_n)$ converges weakly to $(u,v)\in H^1(\T)\times H^1(\T)$.
  By the Sobolev embedding theorem $H^1(\T)$ embeds 
  compactly to the space of H\"older continuous functions with exponent less than $1/2$ 
  and thus we have $v_n\to v$ strongly in $C^{0,\tau}(\T)$ for any $\tau<1/2$.
  Furthermore, it follows that
  \begin{equation}
    \sup_{t\in\T} |Q_nv_n(t) -v_n(t)| \leq
    	\sup_{t\in\T} N \int_{I^N_{j(t)}} |v_n(t)-v_n(t')| dt' \\
    \leq N^{-\tau} \norm{v_n}_{C^{0,\tau}} \\
  \end{equation}
  where $j(t)$ is such that $t\in I^N_{j(t)} = [\frac{j(t)}N,\frac{j(t)+1}{N})$.
  Now we see 
  \[
    \norm{Q_nv_n-v}_{L^2} \leq \norm{Q_nv_n -v_n}_{L^\infty} + \norm{v_n-v}_{L^2} \to 0
  \]
  as $n\to \infty$. The immediate consequence is that
  \[
    \lim_{n\to\infty} \int_\T -N^{1-\alpha} \log(\ep^2+(Q_nv_n)^2) dt
    = 
    -\delta_{\alpha,0} \int_\T \log(\ep^2+v^2) dt
  \]
  when $\alpha \geq 1$.
  Moreover, it also holds that
  \[
    \lim_{n\to\infty} \int_\T \frac{1}{4\ep}(1-v_n)^2 dt 
    = \int_\T \frac{1}{4\ep}(1-v)^2 dt .
  \]
  Let us now consider the condition $(i)$ in Definition \ref{def:gconv} of $\Gamma$-convergence.
  Assume that $(u_n,v_n)\to (u,v)$ in weak topology of $H^1(\T)\times H^1(\T)$ as $n\to\infty$.
  Since $Dv_n\to Dv$ weakly in $L^2(\T)$ and since a norm is lower semicontinuous we deduce that
  \begin{equation}
    \label{disc_ep_dv}
    \int_\T \ep |Dv|^2dt \leq \liminf_{n\to\infty} \int_\T \ep |Dv_n|^2dt.
  \end{equation}
  Without losing any generality we may also assume 
  $\Fda(u_n,v_n) \leq C <\infty$ since otherwise there is nothing
  to prove. Hence in particular $\int_\T |\td u_n|^2 dt \leq C/\ep^2$
  and 
  \[
    \lim_{n\to\infty}\int_\T |v^2-v_n^2||\td u_n|^2 dt 
    \leq \lim_{n\to\infty} \norm{v^2-v_n^2}_\infty \cdot \frac{C}{\ep^2} = 0.
  \]
  Consequently, we obtain
  \begin{eqnarray}
    \label{disc_nonlinear_term}
     \int_\T (\ep^2+v^2)|\td u|^2 dt 
    & \leq & \liminf_{n\to\infty} \int_\T(\ep^2+v_n^2)|\td u_n|^2 dt + \lim_{n\to\infty} \int_\T (v^2-v_n^2)|\td u_n|^2dt \nonumber \\ 
    & = & \liminf_{n\to\infty} \int_\T(\ep^2+v_n^2)|\td u_n|^2 dt.
  \end{eqnarray}
  due to lower semicontinuity of the norm $\norm{\cdot} = \norm{\sqrt{\ep^2+v^2} \;\cdot}_{L^2}$
  and the weak convergence of $\td u_n$.
  By combining all inequalities above it follows that
  $\Fbe(u,v) \leq \liminf_{n\to\infty} \Fda(u_n,v_n)$.
  This proves $(i)$ in Definition \ref{def:gconv}.

  For condition $(ii')$ in Definition \ref{def:gconv} we note that
  for an arbitrary $(u,v)\in H^1(\T)\times H^1(\T)$ there exists a sequence
  $(u_n,v_n)\in PL(n)\times PL(n)$ converging to $(u,v)$ in the strong topology.
  It is easy to see that one can then change $\liminf$ into $\lim$ and inequalities to equalities
  in \eqref{disc_ep_dv} and \eqref{disc_nonlinear_term}. This yields the claim.
\end{proof}


\subsection{Convergence with respect to $\ep$}

Let us prove the $\Gamma$-convergence for a modified functional. 
Define $\Xi_\ep: L^1(\T)\times L^1(\T)\to (-\infty,\infty]$ as
\[
  \Xi_\ep(u,v) =  \left\{
  \begin{array}{ll}
    \Fe^1(u,v) & \textrm{when } (u,v) \in X, \\
    \infty & \textrm{otherwise.}
  \end{array}
  \right. 
\]

\begin{theorem}
\label{Xi_conv}
It holds that $\MS = \Glim_{\ep\to 0} \Xi_\ep$
in the strong topology of $L^1(\T)\times L^1(\T)$.
\end{theorem}

\begin{proof}
First we show that condition $(i)$ of Definition \ref{def:gconv} holds.
Suppose that $\lim_{\ep\to 0} (u_\ep,v_\ep) = (u,v)$ in $L^1(\T)\times L^1(\T)$.
As in the previous $\Gamma$-convergence proofs we may assume without losing any generality that 
\begin{equation}
  \label{aux_xi_gconv}
  \liminf_{\ep\to 0} \Xi_\ep(u_\ep,v_\ep) \leq C < \infty
\end{equation}
and $(u_\ep,v_\ep)\in X$. 
By using the same technique as in 
Lemma \ref{infFep} we can show a lower bound
\begin{equation}
  \label{aux_xi_gconv2}
  \Xi_\ep(u_\ep,v_\ep) \geq \pone(v_\ep) + \int_\T \frac{1}{8\ep} (1-v_\ep)^2dt \geq -C'\ep
\end{equation}
for some constant $C'>0$. Moreover, inequalities \eqref{aux_xi_gconv} and \eqref{aux_xi_gconv2} yield
\begin{equation}
  \label{convto1}
  \int_\T \frac{1}{8\ep} (1-v_\ep)^2dt \leq C+C'\ep
\end{equation}
and hence
in particular $v_\ep \to \onef$ in $L^2(\T)$ as $\ep\to 0$ and $v=1$ a.e. Since $0\leq v_\ep\leq 1$ 
we have by Lemma \ref{logconv} that
$\lim_{\ep\to 0} \pone(v_\ep) = 0$. 

Let us next show that 
$\lim_{\ep\to 0} \ptwo(u_\ep,v_\ep) = 0$. Again since $v_\ep\leq 1$ we have
\begin{equation}
  \label{aux_xi_gconv3}
  |\ptwo(u_\ep,v_\ep)| \leq C \ep^q b_\ep \int_\T |\td u_\ep|dt 
  \leq C \ep^q b_\ep \sqrt{\int_\T |\td u_\ep|^2dt}
\end{equation}
where $b_\ep = \abs{\int_\T u_\ep dt}$.
Since $u_\ep\to u$ in $L^1(\T)$ then also $b_\ep\to |\int_\T u dt|< \infty$.
By assumption it holds that $\ep^2 \int_\T |\td u_\ep|^2dt \leq C$ so that since $q>1$ we obtain
that the right-hand side of inequality \eqref{aux_xi_gconv3} converges to zero as $\ep\to 0$.

Now Theorem \ref{attheorem} implies
\[   
  \MS(u,v) \leq
  \liminf_{\ep\to 0} \Fat(u_\ep,v_\ep) + 
  \lim_{\ep\to 0} \pone(v_\ep) + \lim_{\ep\to 0} \ptwo(u_\ep,v_\ep) = \liminf_{\ep\to 0} \Xi_\ep(u_\ep,v_\ep).
\]
This yields condition $(i)$. 

Next let us consider condition $(ii)$. By Theorem \ref{attheorem}
for any $(u,v)\in L^1(\T)\times L^1(\T)$ there exists a sequence $\{(u_\ep,v_\ep)\}\in X$ such that
$\limsup_{\ep\to 0} \Fat(u_\ep,v_\ep)\leq \MS(u,v)$.
By assuming that $\MS(u,v)$ is bounded we obtain inequality \eqref{convto1} for $v_\ep$
and Lemma \ref{logconv} yields 
$\lim_{\ep\to 0}\pone(v_\ep)=0$. Since also $\ep^2\int_\T |Du_\ep|^2dt \leq C$
the convergence of $\ptwo(u_\ep,v_\ep)$ to zero follows from the estimate
\[
  |\ptwo(u_\ep,v_\ep)| \leq C\left( \ep^q b_\ep \int_\T |D u_\ep|dt + \ep^{2q} b_\ep^2\right)
\]
where $b_\ep = |\int_\T u_\ep dt|$.
Finally we can conclude that
\[  
  \limsup_{\ep\to 0} \Xi_\ep(u_\ep,v_\ep)
  = \limsup_{\ep\to 0} \Fat(u_\ep,v_\ep) + \lim_{\ep\to 0} \pone(v_\ep)
  +\lim_{\ep\to 0} \ptwo(u_\ep,v_\ep) \leq MS(u,v).
\]
This proves $(ii)$ in Definition \ref{def:gconv} and hence the claim follows.
\end{proof}

\begin{theorem}
  \label{fbegconv}
  We have that $\MS = \Glim_{\ep\to 0} \Fbe$ in $L^1(\T)\times L^1(\T)$ for any $\alpha\geq 1$.
\end{theorem}

\begin{proof}
  We prove only the case when $\alpha=1$. For $\alpha>1$ the proof is obtained by leaving out the considerations
  regarding the term $\pone$.
  
  Notice that the functionals $\Fe^1$ and $\Xi_\ep$ differ only in set $X_\ep \setminus X$.
  Obviously since $X \subset X_\ep$ the condition $(ii)$ in Definition \ref{def:gconv} follows immediately from inequality
  $\Fe^1 \leq \Xi_\ep$ and Theorem \ref{Xi_conv}. 

  Let us then consider the condition $(i)$ and let $(u_\ep,v_\ep)\in X_\ep$
  be a sequence converging to $(u,v) \in L^1(\T)\times L^1(\T)$ as $\ep\to 0$. By assuming that sequence
  $\Fe^1(u_\ep,v_\ep)$ is bounded we notice as in the proof of Theorem \ref{Xi_conv} that $v=1$ almost everywhere.
  Consider the folding operation $\Psi_1$ defined in equation \eqref{psifunction}.
  One can easily show that since $0\leq v_\ep \leq 1+30\ep$ 
  we must have $\tilde v_\ep=\Psi_1\circ v_\ep \to v$ in $L^1(\T)$.
  Furthermore, due to Lemma \ref{psilemma} we have 
  \[
    \Xi_\ep(u_\ep,\Psi_1\circ v_\ep) \leq \Fe^1(u_\ep,v_\ep) + \pone(\Psi_1 \circ v_\ep) - \pone(v_\ep)
  \]
  Clearly we have $\pone(\Psi_1 \circ v_\ep) - \pone(v_\ep) \to 0$ as $\ep\to 0$ and thus
  \[
    \MS(u,v) \leq \liminf_{\ep\to 0} \Xi(u_\ep,\tilde v_\ep) \leq \liminf_{\ep\to 0} \Fe^1(u_\ep,v_\ep)
  \]
  which yields the result.
\end{proof}

\begin{lemma}
  \label{L2limsup}
  Let $\alpha\geq 1$ and $u\in \sbv(\T)$.
  For any sequence $\{\ep_j\}_{j=1}^\infty$, $\ep_j\to 0$,
  there exists functions $\{(u_{\ep_j},v_{\ep_j})\}_{j=1}^\infty \subset X$
  converging to $(u,\onef)$ in $L^2(\T)\times L^2(\T)$
  such that
  \begin{equation}
    MS(u,\onef) = \lim_{j\to\infty} \Fbe(u_{\ep_j},v_{\ep_j}).
  \end{equation}
\end{lemma}

\begin{proof}
  Since $u\in\sbv(\T) \subset L^\infty(\T)$ we may apply Theorem \ref{attheorem}
  for $p=2$ to see that there are $\{(u_{\ep_j},v_{\ep_j})\}_{j=1}^\infty \subset X$ 
  converging to $(u,\onef)$ in $L^2(\T)\times L^2(\T)$ such that
  $MS(u,\onef) = \lim_{j\to\infty} \Fat(u_{\ep_j},v_{\ep_j})$.  Following the proof of
  Theorem \ref{Xi_conv} we can show
  $\lim_{j\to\infty} (L_{\ep_j}(v_{\ep_j}) + S^q_{\ep_j}(u_{\ep_j},v_{\ep_j})) = 0$ and the claim follows.
\end{proof}


\subsection{Convergence of minimas}
\label{sec:equicoer}

Let us show the equi-coercivity for the sequences of functionals studied above.

\begin{lemma}
  \label{aux_equicoer}
  Assume that $v\in H^1(a,b)$, $a<b$, $a,b\in \R$ and
  $\max_{t\in [a,b]} v(t) - \min_{t\in [a,b]} v(t) \geq T$.
  Then it holds that
  $\int_a^b |Dv(t)|^2dt \geq \frac{T^2}{b-a}$.  
\end{lemma}

\begin{proof}
By the Sobolev embedding theorem $v$ can be extended to a continuous function on $[a,b]$.
Denote by $t_+$ and $t_-$ points in $[a,b]$ where
\[
  v(t_+) = \max_{t\in [a,b]} v(t) \quad {\rm and} \quad v(t_-) = \min_{t\in [a,b]} v(t).
\]
Without losing the generality we may assume that $t_+>t_-$. Then 
using the fundamental theorem of calculus we see that
\[
  T \leq v(t_+) -v(t_-) \leq \int_{t_-}^{t_+} |Dv(t)| dt 
  \leq \norm{Dv}_{L^2(a,b)} \sqrt{b-a}.
\]
This proves the statement.
\end{proof}

Let us next prove the equi-coerciveness of the functionals $\Fbe$.
\begin{theorem}
  \label{equicoer_ep}
  Let $\alpha\geq 1$, $C>0$, and 
  $(u_\ep,v_\ep)\in H^1(\T)\times H^1(\T)$, $\ep>0$, be a sequence such that 
  \begin{equation}
    \label{f_unif_bded}
    \Fbe(u_\ep,v_\ep) + \norm{Au_\ep-m}_{L^2}^2 \leq C.
  \end{equation}
  Then there exists subsequence 
  $\{(u_{\ep_j},v_{\ep_j})\}_{j=1}^\infty$ which converges
  in $L^1(\T)\times L^1(\T)$.
\end{theorem}

\begin{proof}
The proof is principally the same for both cases $\alpha=1$ and $\alpha>1$.
First notice that the assumption \eqref{f_unif_bded} yields $\norm{v_\ep-\onef}^2_{L^2} \leq C\ep$
and hence the convergence of $v_\ep$ to $\onef$ in $L^1(\T)$ is clear. 
The case for $u_\ep$ follows by considering carefully how the convergence of $v_\ep$
takes place.
Let us first fix some $\ep > 0$ and divide the domain $\T$ into $K = \lfloor\frac 1\ep\rfloor+1$
half-open intervals $I^K_k = [\frac k K, \frac{k+1}{K})$ where $\lfloor\frac 1\ep\rfloor$ denotes the largest
integer less or equal to $1/\ep$. Moreover, denote
\[
  {\mathcal I}_\ep = \left\{ 0\leq k < K \; : \; \max_{t\in I^K_k} v_\ep(t) - \min_{t\in I^K_k} v_\ep(t) \geq \frac 14\right\}.
\]
From Lemma \ref{aux_equicoer} and inequality \eqref{f_unif_bded} we deduce that 
\[
  \sharp({\mathcal I}_\ep) \cdot \frac{K}{16} \leq \int_\T |Dv|^2dt \leq \frac{C}{\ep},
\]
and hence $\sharp ({\mathcal I}_\ep) \leq 16 C'$.
Furthermore, denote
\[
  {\mathcal J}_\ep = \left\{k\in\{0,...,N-1\} \setminus {\mathcal I}_\ep \; | \; \min_{t\in I_k^K} v_\ep(t) < \frac{1}{4}\right\}.
\]
If $j\in {\mathcal J}_\ep$ we observe that the minimum of $v_\ep$ on the interval $I^K_j$ is less than $1/4$
but also the oscillation is less than $1/4$. In consequence, we have that 
$I^K_j \subset \{t\in\T \; | \; v_\ep(t) \leq \frac 12\}$
and the boundedness in inequality \eqref{f_unif_bded} yields
\[
  \sharp({\mathcal J}_\ep) \cdot \frac 1 K \leq 4 \int_\T(1-v_\ep)^2dt \leq 16 \ep C
\]
and thus $\sharp({\mathcal J}_\ep) \leq 16 C''$. Consider the union of all intervals
in the complement of ${\mathcal I}_\ep\cup {\mathcal J}_\ep$. Since the number of indeces in ${\mathcal I}_\ep\cup {\mathcal J}_\ep$ is less than $L = \lfloor16(C'+C'')\rfloor + 1$
it follows that its complement $({\mathcal I}_\ep\cup {\mathcal J}_\ep)^c$ can be presented 
as a union of at most $L$ half-open connected intervals $K_j$ so that
\[
  J(\ep) := \bigcup_{0\leq k < K, k \notin {\mathcal I}_\ep\cup {\mathcal J}_\ep} I^K_k
  = \bigcup_{j=1}^L K_j.
\]
Next we obtain $L^q$-boundedness for $u_\ep$ with some $q>1$ by applying
the Poincare inequality on each interval $I^K_j$, $j\in {\mathcal I}_\ep\cup {\mathcal J}_\ep$, and to every $K_j$.
Let $I\subset \T$ be an open connected interval and $p>1$ such that
$\frac{1}{p} = \frac 12 +s$.
By the Poincare inequality we have
\begin{equation}
  \label{ue_poincare}
  \norm{u_\ep-b_\ep}_{L^p(I)} \leq c |I| \norm{Du_\ep}_{L^p(I)}
\end{equation}
where $b_\ep= \frac{1}{|I|} \int_I u_\ep(t)dt$ is the average of $u_\ep$ on $I$ which
by Lemma \ref{dual_bound} satisfies
$\abs{b_\ep} \leq C |I|^{-p}$.
Using triangle inequality to estimate the left-hand side of \eqref{ue_poincare} we obtain
\begin{equation}
  \label{eq:jepjep}
  \norm{u_\ep}_{L^p(I)} \leq c|I| \norm{Du_\ep}_{L^p(I)} + |I|^{\frac 1p} |b_\ep| 
  \leq C(|I| \norm{Du_\ep}_{L^p(I)} + 1).
\end{equation}
Let us now apply the above inequality to
\begin{equation}
  \label{lq_ue_1}
  \norm{u_\ep}_{L^p(\T)} \leq \sum_{k\in {\mathcal I}_\ep\cup {\mathcal J}_\ep} \norm{u_\ep}_{L^p(I^K_k)}+
  \sum_{j=1}^L \norm{u_\ep}_{L^p(K_j)}.
\end{equation}
For any interval $I^K_k$ with $k\in {\mathcal I}_\ep\cup {\mathcal J}_\ep$ we have by the H\"older inequality that
\begin{equation}
  \label{lq_ue_aux1}
  |I^K_k| \norm{Du_\ep}_{L^p(I^K_k)} \leq \ep \norm{Du_\ep}_{L^p(I^K_k)} \leq \ep^s
  \norm{\ep Du_\ep}_{L^2(I^K_k)} \leq C\ep^s. 
\end{equation}
On any interval $K_j$, $1\leq j \leq L$, we know that $v_\ep>\frac 12$. By the boundedness assumption
\eqref{f_unif_bded} we have then
\begin{equation}
  \label{lq_ue_aux2}
  \norm{Du_\ep}_{L^p(\cup_{j=1}^L K_j)} \leq C\norm{\sqrt{\ep^2+v_\ep^2} Du_\ep}_{L^2(\cup_{j=1}^L K_j)} \leq C'.
\end{equation}
Applying inequalities \eqref{eq:jepjep}, \eqref{lq_ue_aux1} and \eqref{lq_ue_aux2} to \eqref{lq_ue_1} yields
\[
  \norm{u_\ep}_{L^p(\T)} \leq C\left( \sum_{k\in {\mathcal I}_\ep\cup {\mathcal J}_\ep} (\ep^s + 1)
  + \sum_{j=1}^L (|K_j|+1) \right)\leq C 
\]
since $L$ and $\sharp ({\mathcal I}_\ep\cup {\mathcal J}_\ep)$ were bounded. 
Notice that this bound is uniform with respect to $\ep$. By the Banach-Alaoglu theorem
we can extract a subsequence, denoted also by $u_\ep$, which converges to some $u\in L^p(\T)$ weakly.

Next, by inequality \eqref{lq_ue_aux2} we obtain $\norm{u_\ep}_{W^{1,q}(J(\ep))} < C$ for some constant $C>0$ independent of $\ep$.
Let us then extract a subsequence $\{u_{\ep_j}\}_{j=1}^\infty$ in the following way: denote $Z_\ep = \{ k/K :  k\in {\mathcal I}_\ep\cup {\mathcal J}_\ep\}\subset \T$.
Let $\{\ep_j\}_{j=1}^\infty$ be such that the set $Z_{\ep_j}$ converges in the Hausdorff topology to some discrete set $Z$ such that $\sharp(Z)\leq L$.
Note that $J(\ep_j)^c$ is included in an $\ep_j$-neighbourhood of the set $Z_{\ep_j}$.
Then it follows that for every $\ell\in \Z_+$ 
and $(1/\ell)$-neighborhood ${\mathcal U}_\ell$ of $Z$
we have $\norm{u_{\ep_j}}_{W^{1,q}({\mathcal U}^c_\ell)} < C$.
Furthermore, by the Banach-Alaoglu theorem we can
extract another subsequence, denoted also by $u_{\ep_{j}}$,
which converges weakly in $W^{1,q}({\mathcal U}^c_\ell)$ for all $\ell$. 
The Sobolev embedding theorem
then yields that this subsequence also converges
strongly in $L^1({\mathcal U}^c_\ell)$. 
We conclude that
\begin{eqnarray*}
  \lim_{j\to \infty} \norm{u-u_{\ep_{j}}}_{L^1(\T)} & \leq & 
  \lim_{j\to \infty} \left(\norm{u-u_{\ep_{j}}}_{L^1({\mathcal U}_\ell)} + 
  \norm{u-u_{\ep_{j}}}_{L^1({\mathcal U}^c_\ell)}\right) \\
  &\leq & \lim_{j\to \infty} \norm{u-u_{\ep_{j}}}_{L^p({\mathcal U}_\ell)} \left(\frac{2L}{\ell}\right)^{1/p'} 
  \leq C \left(\frac{2L}{\ell}\right)^{1/p'}
\end{eqnarray*}
for any $\ell\in \Z_+$ where $1/p+1/p'=1$. Finally, the result follows 
since $\ell$ was arbitrary.
\end{proof}

\begin{theorem}
  \label{equicoer_n}
  Consider fixed $\ep>0$.
  Let $(u_n,v_n)\in PL(n)\times PL(n)$ be a sequence such that 
  \begin{equation}
    \label{f_unif_bded2}
    \Fda(u_n,v_n) + \norm{A_nu_n-m}_{L^2}^2 \leq C 
  \end{equation}
  for some constant $C<\infty$.
  Then there exists subsequence $(u_{n_j},v_{n_j})$ which converges
  weakly in $H^1(\T)\times H^1(\T)$.
\end{theorem}

\begin{proof}
  We see that  $\int_\T \ep^2 |\td u_n|^2dt$  and thus $\norm{u_n}_{H^1}$ are uniformly bounded
  for all $n$. Furthermore, boundedness of 
  $\int_\T \left(\ep |Dv|^2 + \frac{1}{4\ep}(1-v)^2\right)dt$
  yields a bound for $\norm{v_n}_{H^1}$. Then the claim follows by the Banach-Alaoglu theorem.
\end{proof}

\vspace{0.2cm}
{\em Proof to Theorem \ref{pos_result}}.
The claim $(i)$ follows easily by Theorems \ref{gconv_thm1}, \ref{gconv_for_n} and \ref{equicoer_n}.
For the claim $(ii)$ consider the case $\alpha=1$ and let $(\bar u_\ep,\bar v_\ep) \in X_\ep$ be a minimizer
of problem \eqref{eq:fbe+res}. By the equicoercivity theorem \ref{equicoer_ep} we have a subsequence
$\{\ep_j\}_{j=1}^\infty$ such that $\bar v_{\ep_j}$ converges to $\onef$ in $L^1(\T)$ and $\bar u_{\ep_j}$ to some $\bar u$ in $L^1(\T)$.
Let us then prove that $\bar u$ is a minimizer of problem \eqref{eq:ms+res}.
Let $u\in L^1(\T)$ be such that $MS(u,\onef) + R(u)<\infty$. Then clearly $u\in \sbv(\T)$
and by Lemma \ref{L2limsup} there exists a sequence $(u_{\ep_j},v_{\ep_j})_{j=1}^\infty \subset X_\ep$ such that $\lim_{j\to\infty} (u_{\ep_j},v_{\ep_j}) = (u,\onef)$ in $L^2(\T) \times L^2(\T)$ and
$MS(u,\onef) = \lim_{j\to\infty} F_{\ep_j}^1(u_{\ep_j},v_{\ep_j})$. Since the residual $R$ is continuous in $L^2(\T)$ and lower semicontinuous in $L^1(\T)$ we obtain
using Theorem \ref{fbegconv} that
\begin{eqnarray}
  MS(\bar u,\onef) + R(\bar u) & \leq & 
  \liminf_{j\to\infty} F_{\ep_j}^1(\bar u_{\ep_j},\bar v_{\ep_j}) +\liminf_{j\to\infty} R(\bar u_{\ep_j}) \nonumber \\
  & \leq & \liminf_{j\to\infty} \left(F_{\ep_j}^1(u_{\ep_j},v_{\ep_j}) +R(u_{\ep_j})\right) \nonumber \\
  & = & MS(u,\onef) + R(u). \label{minimum_found}
\end{eqnarray}
This proves that $(\bar u,\onef)$ is a minimizer of $MS + R$ since $u\in L^1(\T)$ was arbitrary
 function for which \eqref{minimum_found} is finite.
The case $\alpha>1$ follows identically.
\hfill $\Box$

Finally, before the numerical results, we present two remarks on issues discussed in the introduction.

\begin{remark}
\label{k_not_equal_a}
Let us consider the case when $\alpha\neq \kappa$ in equation \eqref{intro_noise}.
Then the minimization problem of finding MAP estimates can be written as
  \begin{eqnarray*}
    \min_{u,v\in PL(n)} \int_\T &\bigg( &-N^{1-\alpha}\log(\ep^2+(Q_nv)^2) + (\ep^2+(Q_nv)^2)|\td u|^2 \\
   & & +\ep |Dv|^2 + \frac{1}{4\ep} (1-v)^2 + N^{\kappa-\alpha}|P_n A u - m_n|^2 \bigg) dt.
  \end{eqnarray*}
 In consequence, the residual term becomes over or under weighted in the limit regardless of the particular choices  of $\alpha$ or $\kappa$.

If $\alpha = \kappa\geq 1$ then our results are summarized in Theorem \ref{pos_result}
and the MAP estimates converge.
\end{remark}

\begin{remark}
\label{remarkb}
Let us consider statement (b') in the introduction.
When $\alpha\geq 1$,
\begin{eqnarray*}
  \expec \left(\bra V^\alpha_{n,\ep}-\onef,\phi\cet^2_{L^2} \right)
  & = & N^{-\alpha} \norm{\left(\frac{1}{4\ep}I-\ep \Delta\right)^{-1/2} \phi}_{L^2}^2 
  \quad {\rm and} \\ 
  \expec \left(\bra U^\alpha_{n,\ep},\phi\cet^2_{L^2} \right)
  & \leq & \ep^{-2} N^{-\alpha} \norm{\td^{-1} \phi}_{L^2}^2
\end{eqnarray*}
where $\phi\in L^2(\T)$.
Applying this for the orthonormal basis $\{e_j\}_{j=-\infty}^\infty\subset L^2(\T)$, $e_j(t) = \exp(2\pi ijt)$,
yields easily with a fixed $\ep>0$ that
\[
  \lim_{n\to\infty} \expec \norm{V^{\alpha}_{n,\ep}-\onef}_{L^2}^2 = 0 \quad {\rm and} \quad
  \lim_{n\to\infty} \expec \norm{U^\alpha_{n,\ep}}_{L^2}^2 = 0.
\]
In particular, this implies that the random variable $(U^\alpha_{n,\ep},V^{\alpha}_{n,\ep})$
in $L^2(\T)\times L^2(\T)$ 
converges to $({\bf 0},\onef)$ in probability as $n\to\infty$ \cite{Kallenberg}.
\end{remark}

\section{Numerical considerations}
\label{sec:example}

In this section we study the qualitative behaviour of MAP estimates by giving a
numerical example with the scaling $\alpha = 1$. 
Our purpose is to demonstrate that the MAP estimates 
do behave numerically in a similar manner 
in all discretizations, i.e., for different choices of parameter $n$.
This can be expected given the results in Theorem \ref{pos_result}.

The numerical simulations for convergence of the CM estimates are demonstrated in \cite{tapio} in the case $\alpha=0$.

\subsection{The model problem}

We consider a Bayesian deblurring problem with linear operator 
$A = (I-\Delta)^{-s/2} : L^2(\T) \to H^s(\T)$
for a given $0<s<1/2$. Notice that this operator satisfies condition $\norm{u}_{H^{-s}} \leq \norm{Au}_{L^2}$ for any
$u\in H^{-s}(\T)$. We assume the measurements to be obtained via projections
$P_n f = \sum_{j=-N}^N \bra f,e_j\cet_{H^{-1}\times H^1} e_j$
for any $f\in H^{-1}(\T)$ where the $L^2$-orthonormal basis functions $\{e_j\}_{j=-N}^N$ are
$e_j(t) = \exp\left(-2\pi i jt\right)$
for $t\in [0,1)$. It is straight-forward to show that projections $P_n$ are proper measurement projections in the sense
of Definition \ref{def:pk}.

Let us now introduce some notation. For any $n\in\N$ we denote $\phi_j^n\in PL(n)$ a function such that $\phi_j^n(k/n)=\delta_{kj}$.
The basis $\{\phi_j^n\}_{j=1}^n\subset PL(n)$ is called the {\em roof-top functions}.
Let then ${\bf B}_n\in \R^{N\times N}$ be a matrix such that $({\bf B}_n)_{jk} = \bra \phi_j^n,\phi_k^n\cet_{L^2}$
for $1\leq j,k\leq N$.

In the following we use bolded symbols to denote the coefficients of any function $f\in PL(n)$ presented in the roof-top basis, i.e.,
if $f = \sum_{j=1}^N {\bf f}_j \phi_j^n$ then ${\bf f} = ({\bf f}_1, ... , {\bf f}_N)\in \R^N$. Furthermore, denote by ${\bf D}_n$
and ${\bf Q}$ the matrix presentations of operators $\td, Q:PL(n) \to PC(n)$, respectively, from the roof-top basis to the basis of piecewise constant functions
$\{\chi_{I^N_j}\}_{j=1}^N\subset PC(n)$ where $I^N_j = [j/N,(j+1)/N)$. Furthermore, denote the matrix
$\Lambda(\bov) = {\rm diag}(\ep^2+({\bf Q}\bov)^2_j) \in \R^{N\times N}.$
With these notations the functional $\Fda$ written in terms of the coeffients in the roof-top basis functions has the form
\begin{eqnarray}
  \label{disc_prob}
   {\bf F}^\alpha_{\ep,n}(\bou,\bov) & = & -N^{-\alpha}\log(\det \Lambda(\bov)) + \frac 1 N({\bf D}_n \bou)^T \Lambda(\bov) ({\bf D}_n \bou) \nonumber \\
  & & + \frac{\ep}{N} \norm{{\bf D}_n \bov}^2_2 + \frac{1}{4\ep}(\onef-\bov)^T {\bf B}_n (\onef-\bov) + \frac 1{\sigma^2} \norm{{\bf A}_n \bou - \bom}_2^2 
\end{eqnarray}
where ${\bf A}_n\in \R^{(2N+1)\times N}$ maps $\bou$ to the coefficients of $P_nAu$ in the basis $\{e_j\}_{j=-N}^N$. The components of matrix ${\bf A}_n$ satisfy
$({\bf A}_n)_{jk} = \bra e_j, (I-\Delta)^{-s/2} \phi^n_j \cet_{L^2}$
for $-N\leq j\leq N$ and $1\leq k\leq N$.

\subsection{Computational methods}

Because of the non-quadratic terms in ${\bf F}^\alpha_{\ep,n}$ we have chosen
to implement an alternate minimization scheme (see e.g. \cite{bourdin}). 
The convergence of such a method is studied in \cite{NN}
in a setting without the logarithm term $L_\ep$.
Producing a convergence result in our case lies outside the focus of this section.

Let us now write in pseudo-code how the minimizers are achieved:

\begin{itemize}
\item[(1)] Initialize $\bou^0,\bov^0\in \R^N$ and set $j:=1$.
\item[(2)] Solve the equation
 $\left(\frac 1 N {\bf D}_n^T \Lambda(\bov^{j-1}) {\bf D}_n + {\bf A}_n^T {\bf A}_n^T\right) \bou 
  = {\bf A_n}^T \bom$
and set $\bou^j = \bou$.
\item[(3)] Solve the minimization problem
\begin{eqnarray*}
  \min_{\bov\in \R^N} & &\left( -N^{\alpha-1}\log(\det \Lambda(\bov)) + \frac 1 N({\bf D}_n \bou_j)^T \Lambda(\bov) ({\bf D}_n \bou_j)\right. \\
  & & \left.+ \frac{\ep}{N} \norm{{\bf D}_n \bov}^2_2 + \frac{1}{4\ep}(\onef-\bov)^T {\bf B}_n (\onef-\bov) \right)
\end{eqnarray*}
and set $\bov^j = \bov$.
\item[(4)] If $(\bou^j,\bov^j)$ satisfies
${\bf F}^\alpha_{\ep,n}(\bou^{j},\bov^{j}) \leq {\bf F}^\alpha_{\ep,n}(\bou^{j-1},\bov^{j-1}) -\delta$
go to step (2); else stop.
\end{itemize}

\subsection{Results}

\begin{figure}[h]
\begin{picture}(400,180)(40,10)
\epsfxsize=3.6cm 
\epsfysize=2.5cm
\put(44,100){\epsffile{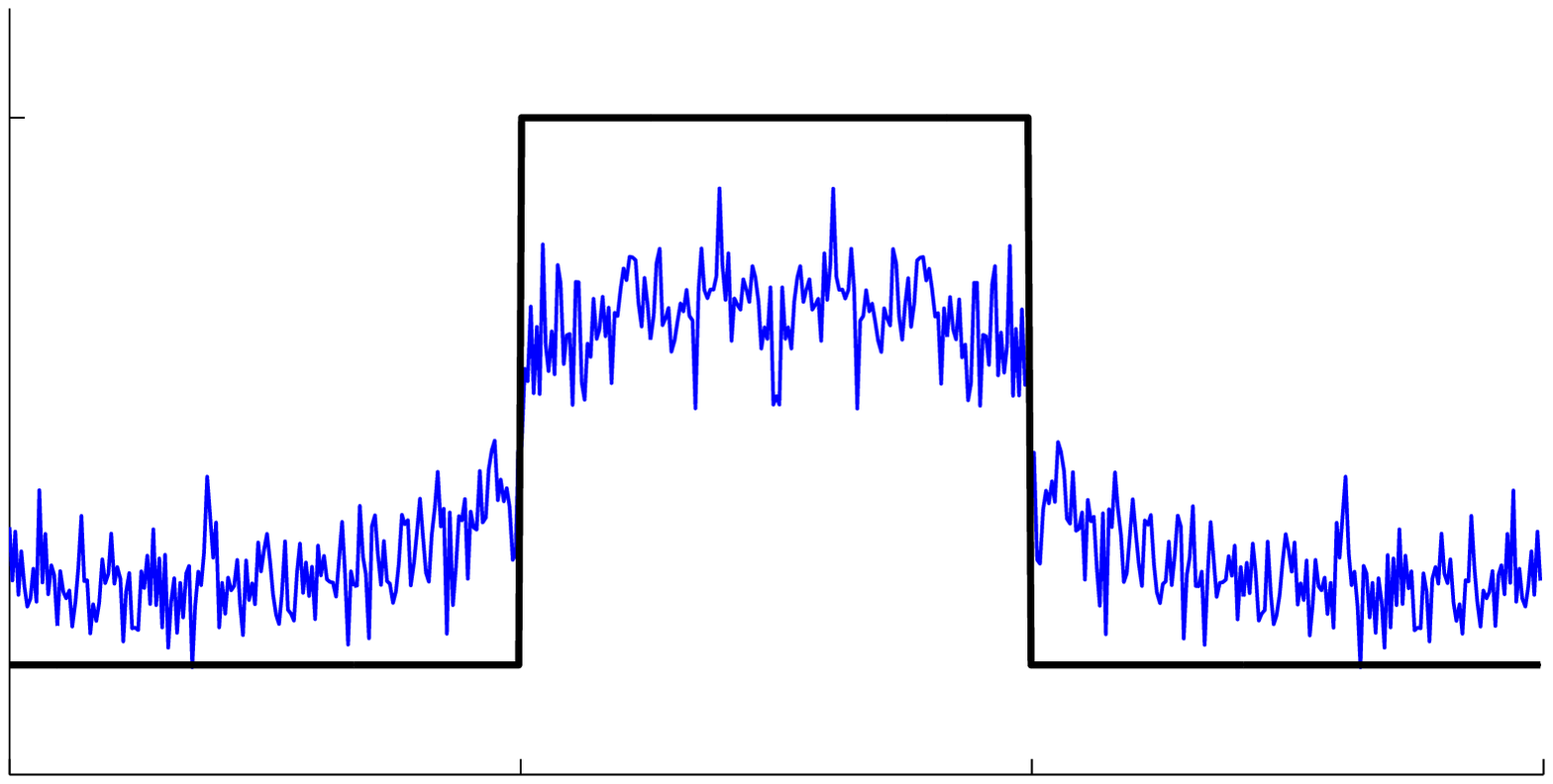}}
\epsfxsize=3.8cm
\put(40,175){$N=512$}
\put(165,100){\epsffile{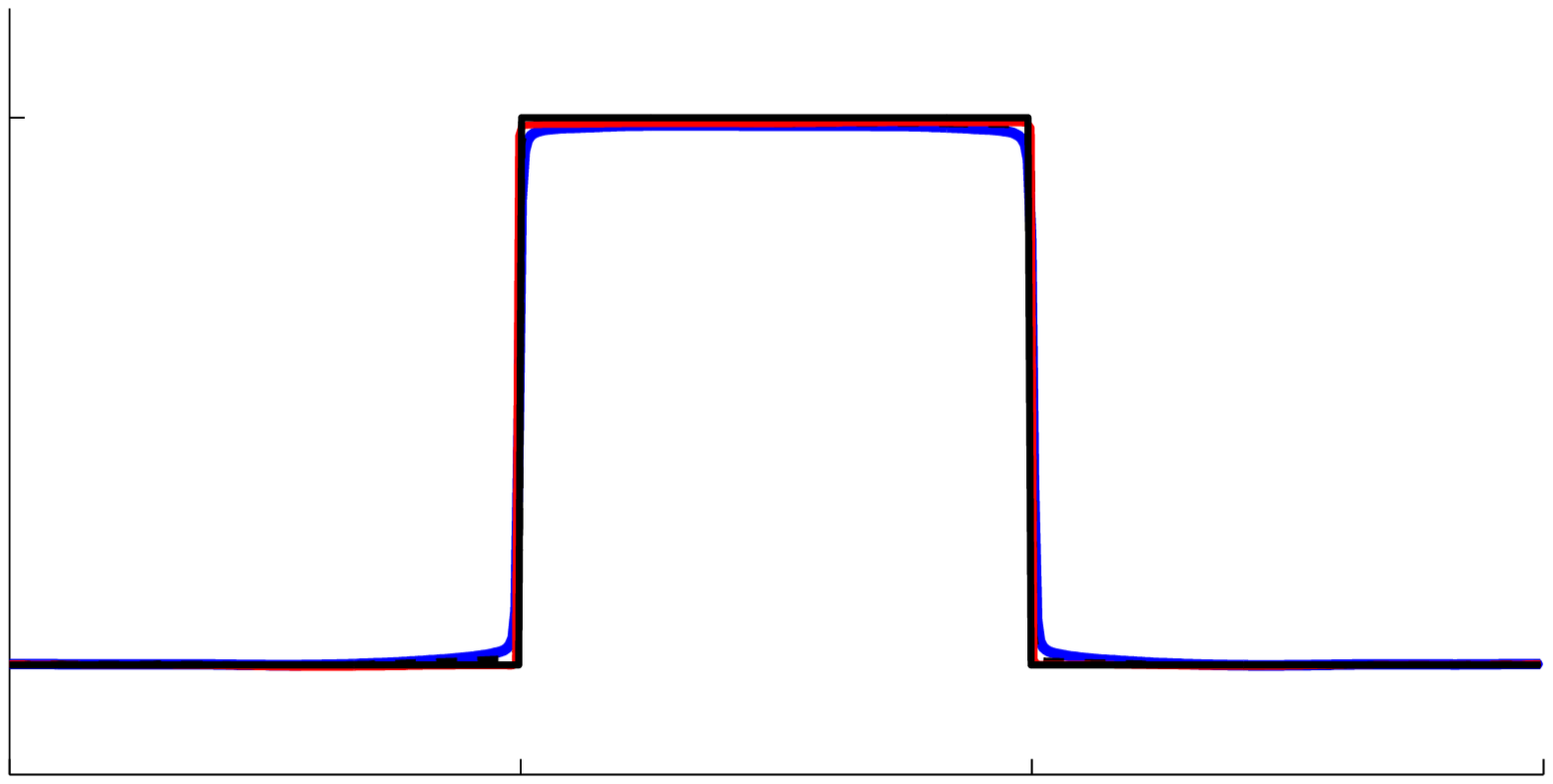}}
\put(290,100){\epsffile{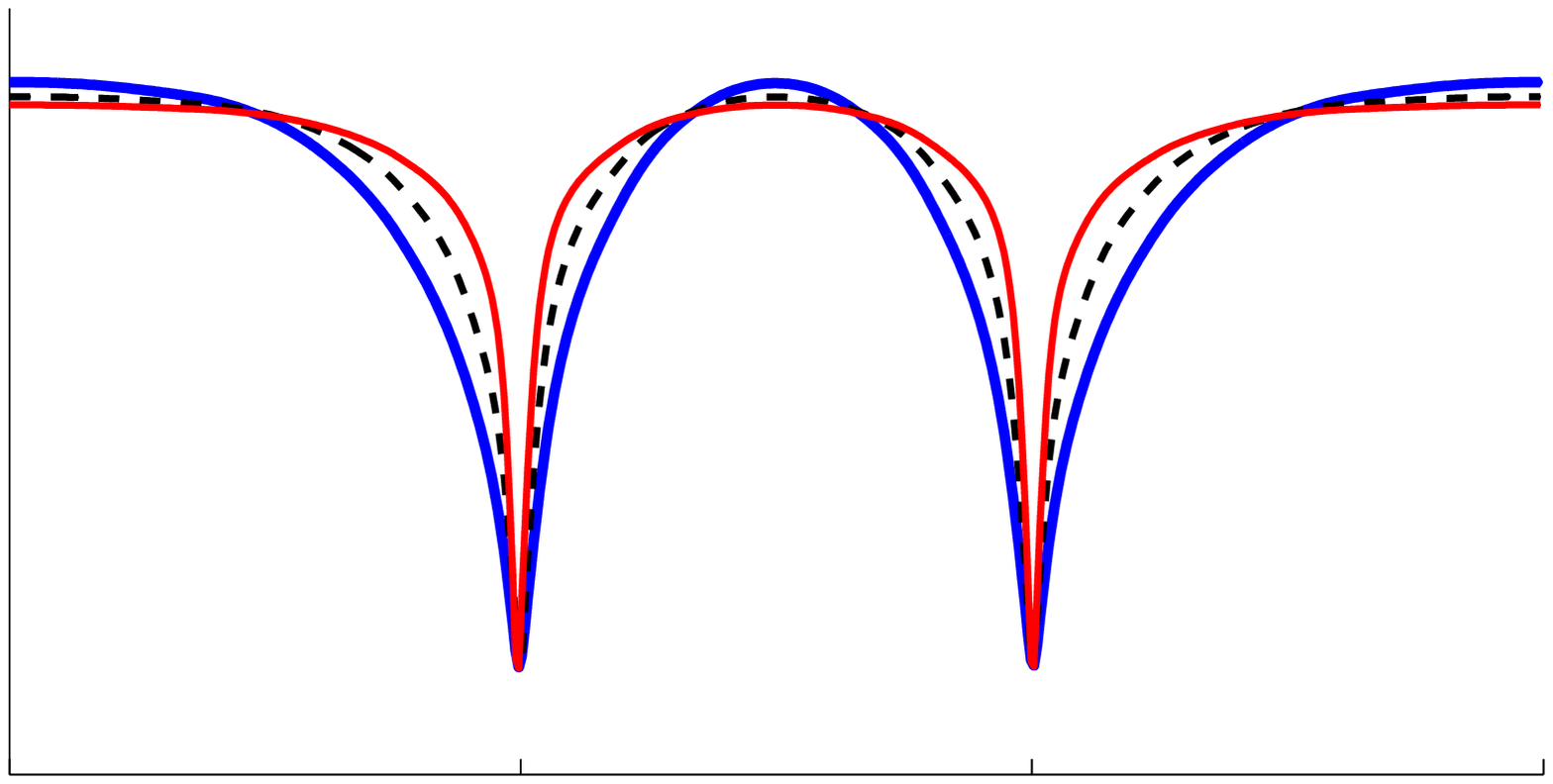}} 

\put(40,10){\epsffile{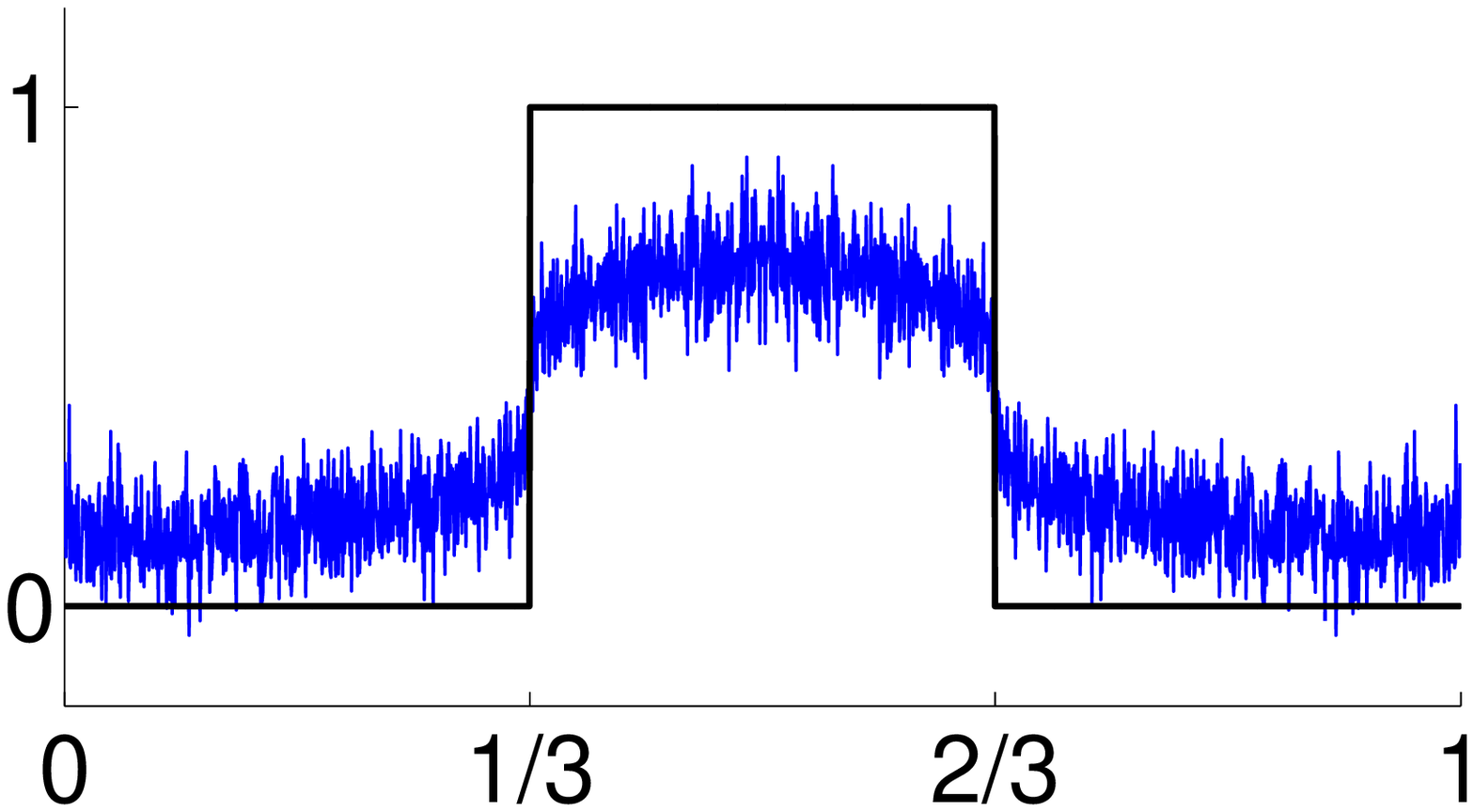}} 
\put(40,85){$N=2048$}
\put(165,10){\epsffile{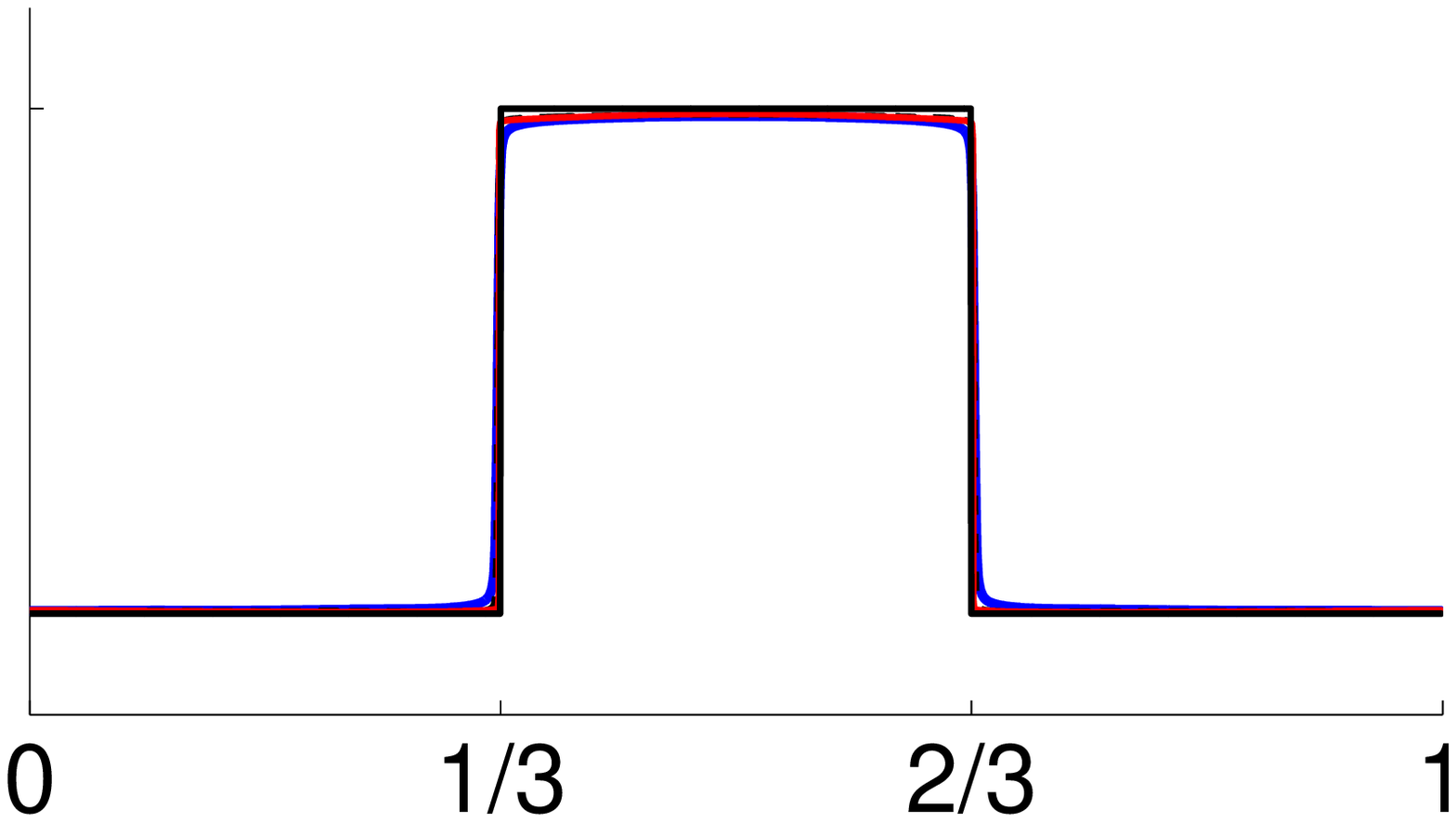}}
\put(290,10){\epsffile{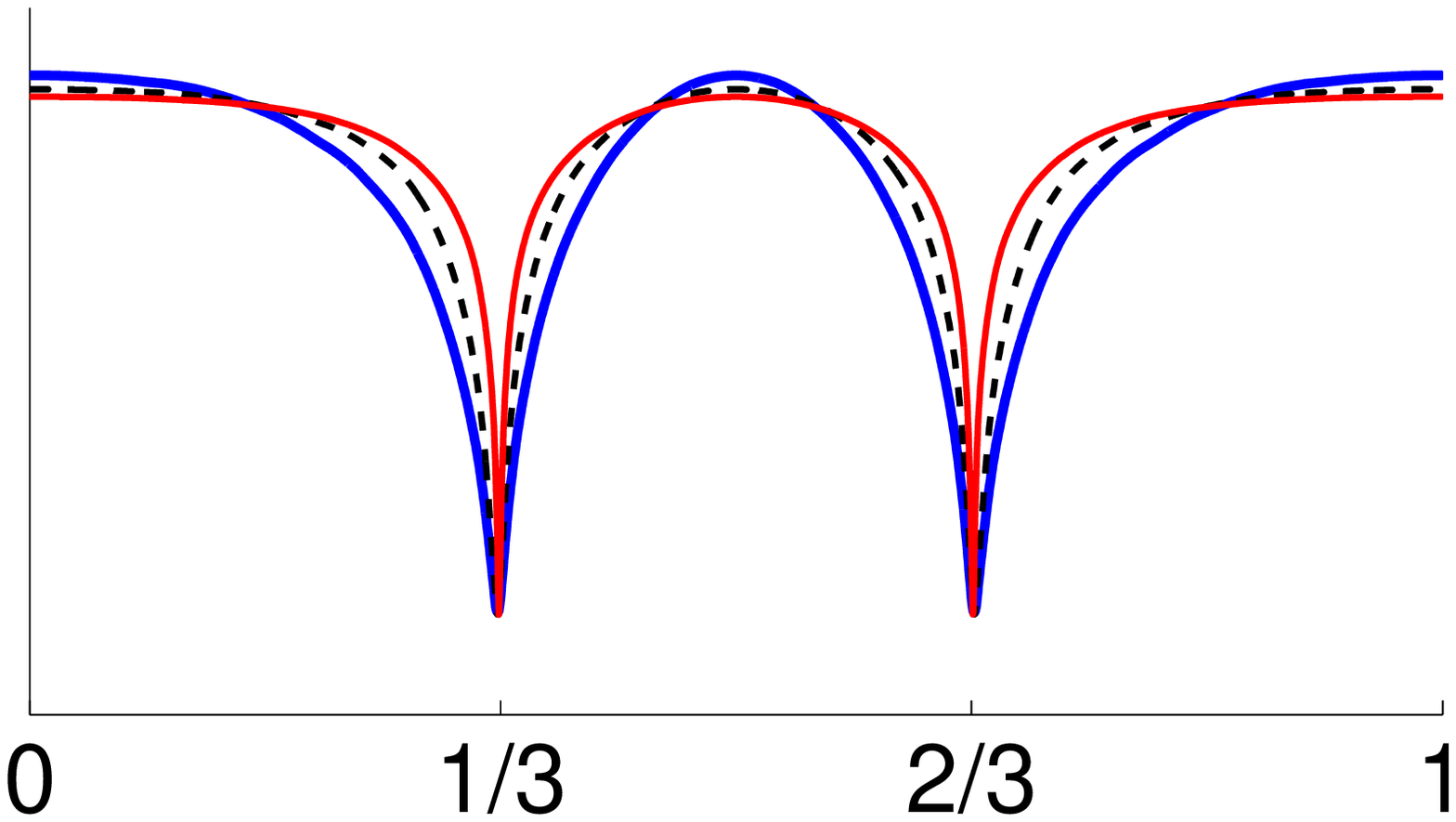}} 

\end{picture}
\caption{A step function.
Left: the true value of $u$ and the noisy measurement $m_n = M_n(\omega_0)$.
Middle: $u_n^{MAP}$ estimates. Right: $v_n^{MAP}$ estimates. The thick, dashed and thin lines represent
reconstruction with sharpness $\ep = 2\times 10^{-2}, 1\times 10^{-2}, 6\times 10^{-3}$, respectively. 
Axis limits are the same in each plot.}
\label{fig_recon}
\end{figure}

\begin{figure}[h]
\begin{picture}(400,180)(40,10)
\epsfxsize=3.6cm 
\epsfysize=2.5cm
\put(44,100){\epsffile{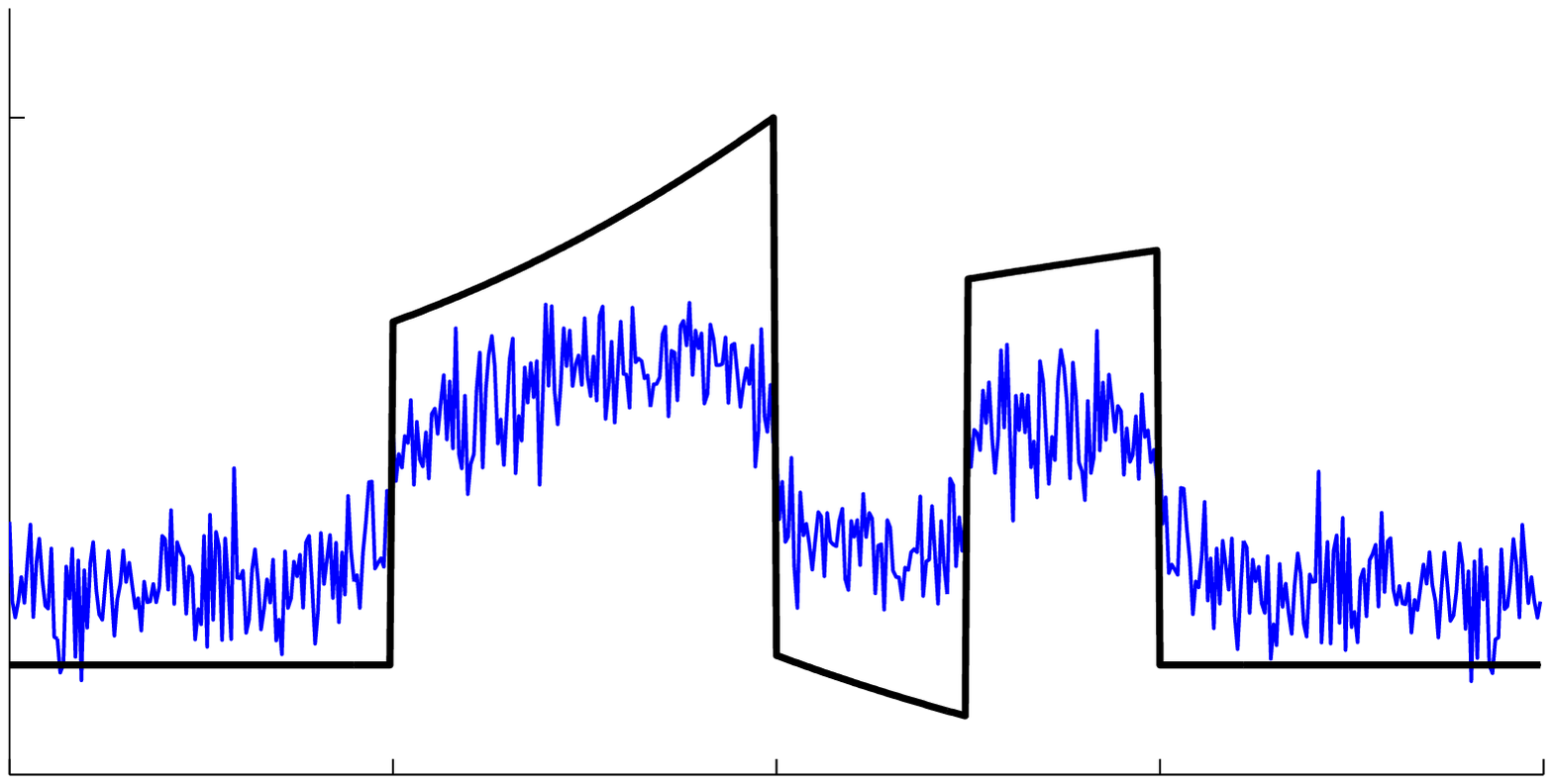}}
\epsfxsize=3.8cm
\put(40,175){$N=512$}
\put(165,100){\epsffile{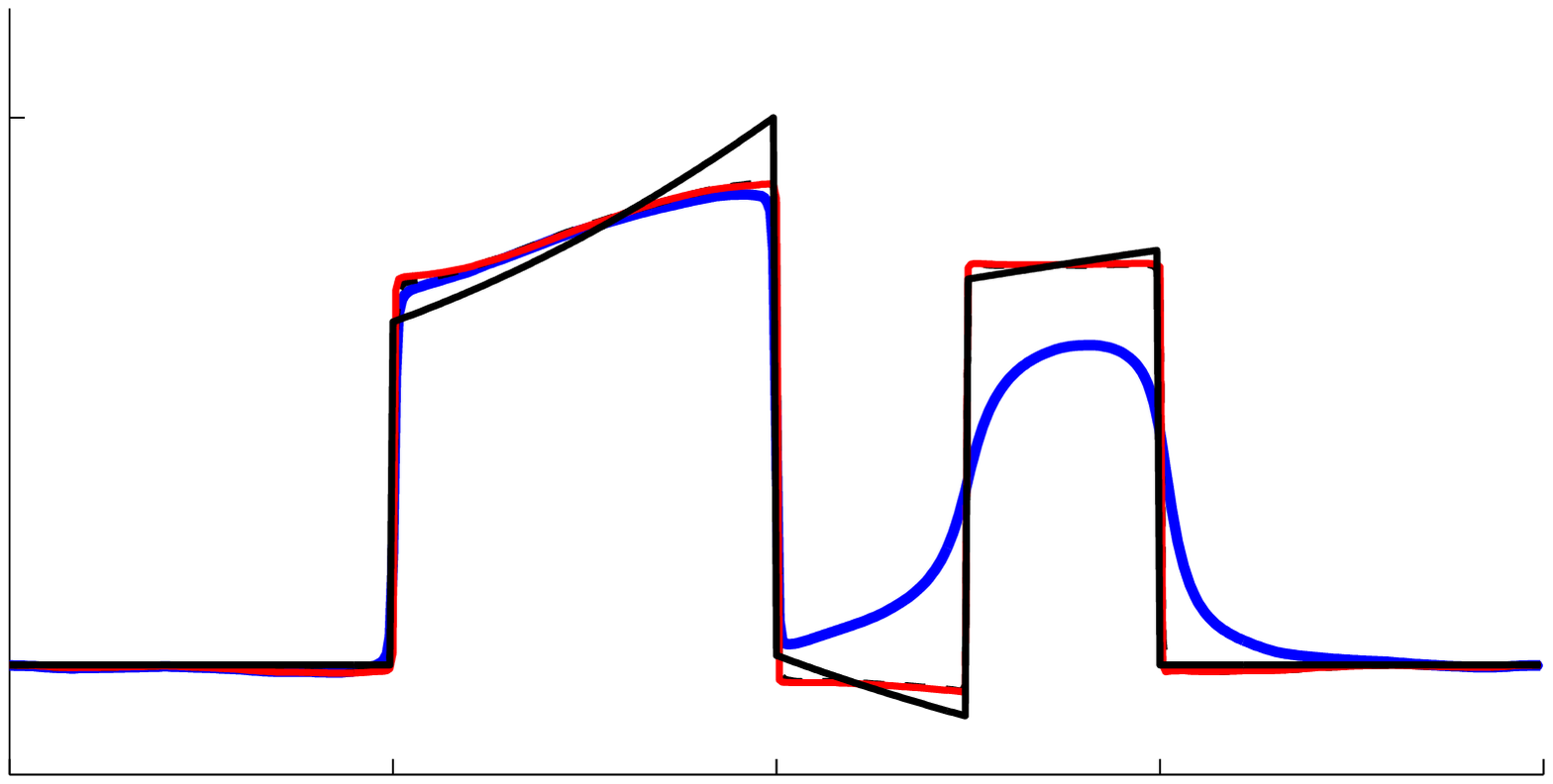}}
\put(290,100){\epsffile{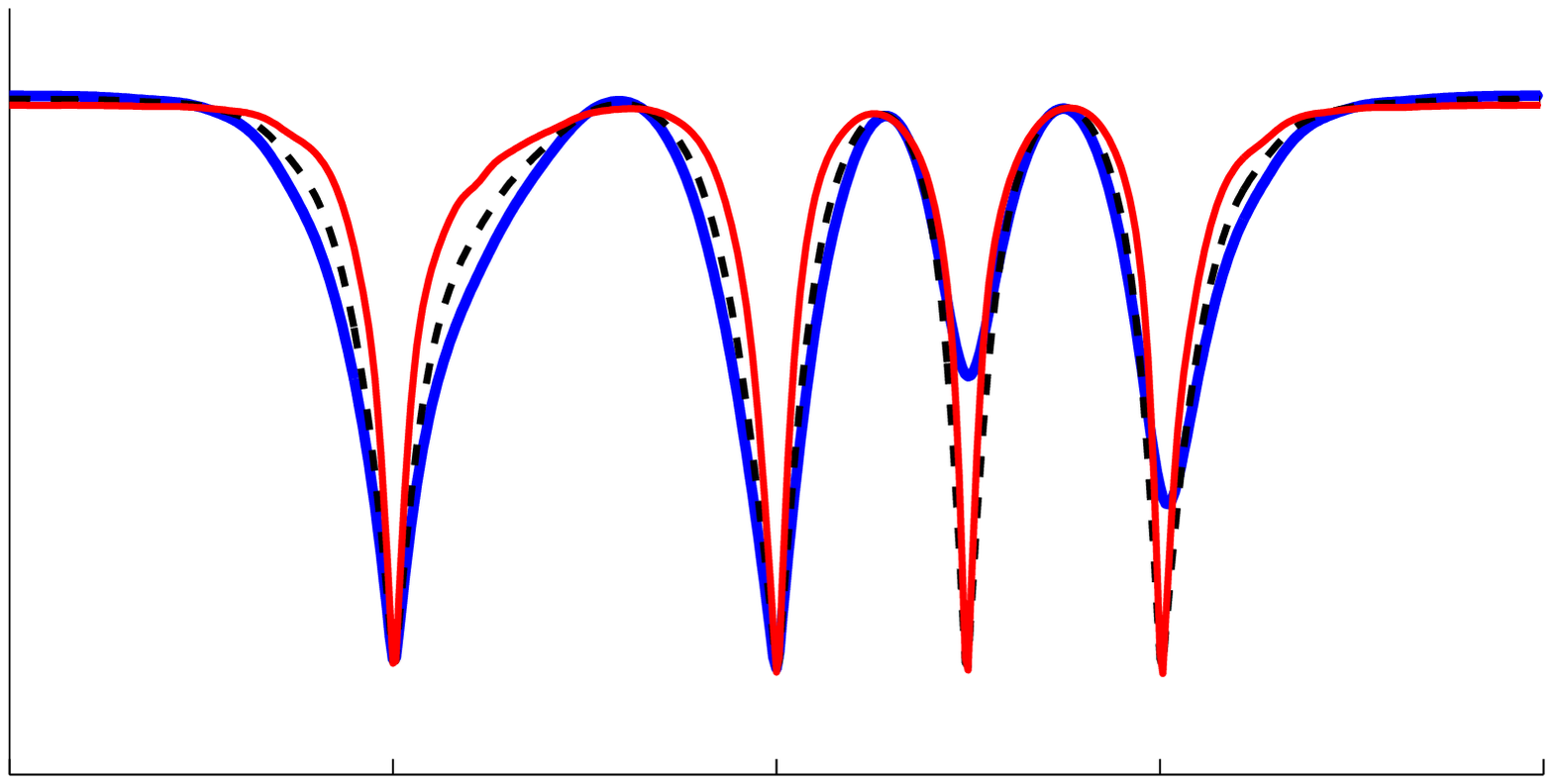}} 

\put(40,10){\epsffile{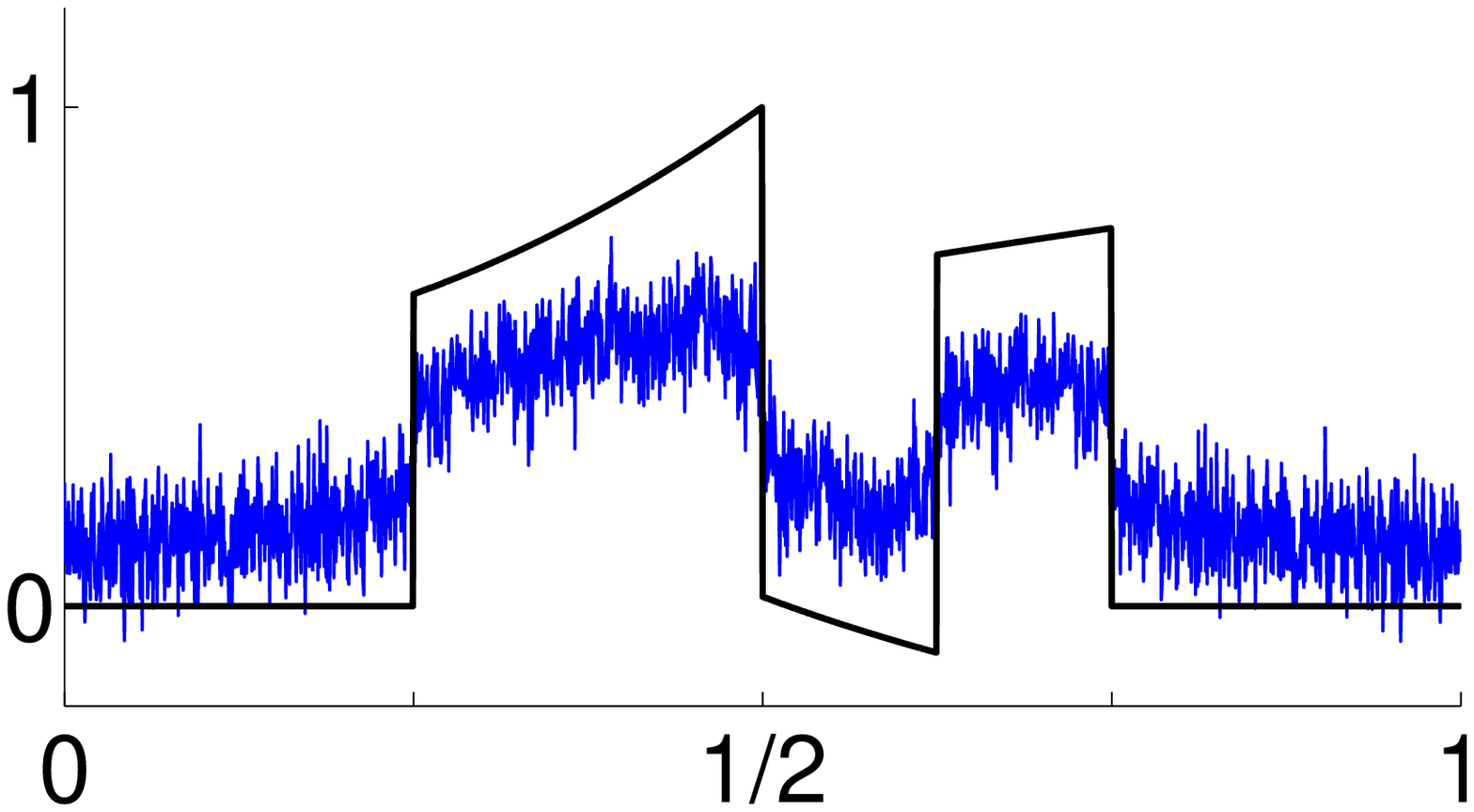}} 
\put(40,85){$N=2048$}
\put(165,10){\epsffile{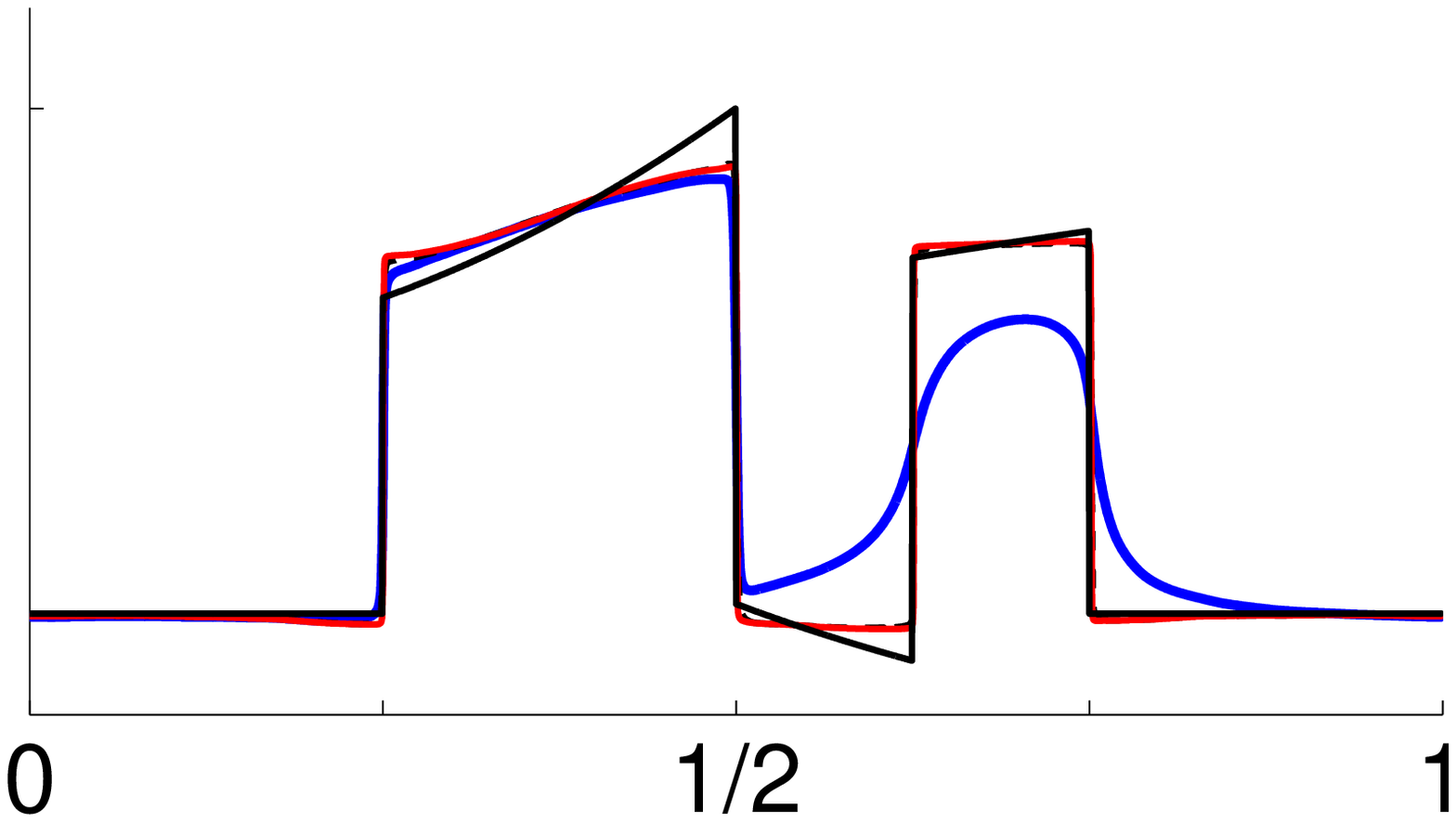}}
\put(290,10){\epsffile{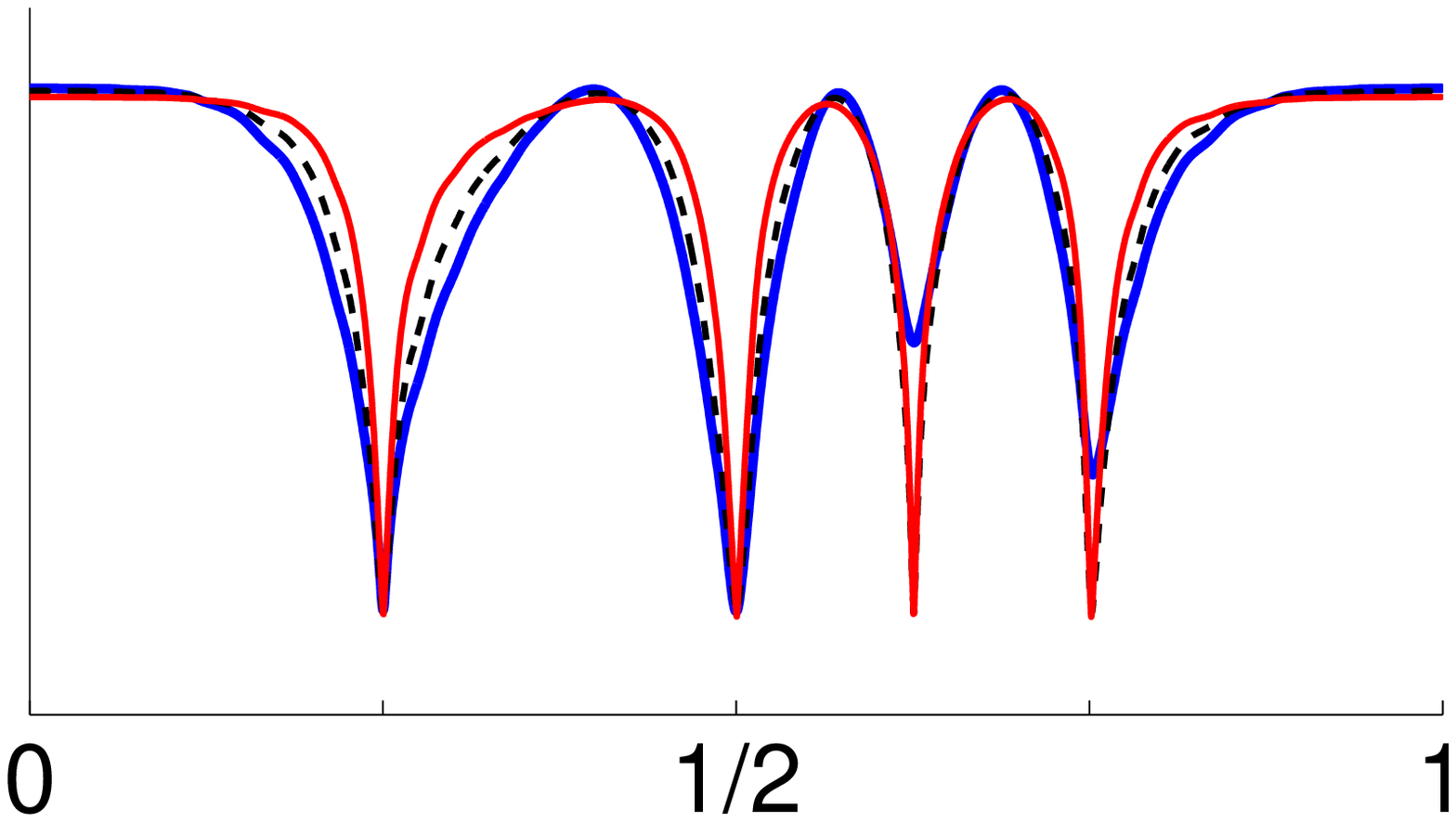}} 

\end{picture}
\caption{A piecewise smooth function.
Left: the true value of $u$ and the noisy measurement $m_n = M_n(\omega_0)$.
Middle: $u_n^{MAP}$ estimates. Right: $v_n^{MAP}$ estimates. The thick, dashed and thin lines represent
reconstruction with sharpness $\ep = 2\times 10^{-2}, 1\times 10^{-2}, 6\times 10^{-3}$, respectively. 
Axis limits are the same in each plot.}
\label{fig_recon2}
\end{figure}

We implemented the problem with operator $A$ having parameter $s=0.35$ and measurement
noise with variance $\sigma = 5 \times 10^{-3}$, i.e., replace $N^{-\kappa}$ in equation \eqref{intro_noise}
with $\sigma N^{-\kappa}$. Furthermore, the scaling of the prior is assumed to be $\alpha=1$.
We used four different sets of data with two true values of $u$
and two discretization sizes $N=512$ and $N=2048$. The MAP estimates were computed with sharpness parameters
$\ep = 2 \times 10^{-2}, 1 \times 10^{-2}, 6 \times 10^{-3}$.
The reconstructions are shown in Figures \ref{fig_recon} and \ref{fig_recon2}.

In Figure \ref{fig_recon} the true value of $u$ is a simple step function. We have weighted the
residual with constant $c = 14$. In Figure \ref{fig_recon2} the true value is piecewise smooth with 
$\sharp(S_u) = 4$ and the residual was weighted with $c=10$.
The initial values in all the computations were vectors $\bou^0 = {\bf 0}$ and $\bov^0 = \onef$.
The step (2) in the algorithm was implemented by using Matlab's {\it backslash} function and
in the step (3) we used a gradient-descent method by choosing step-sizes with a line search algorithm. 
The minimization in step (3) was stopped when either no satisfying step-size was found or the values of the functional
did not change by high accuracy. All computations were stopped at 50 iteration.

We perform all the computations with Matlab 7.6 running in a desktop PC computer with Dual Intel Xeon
processor running at 2,80 GHz and 4 GB of RAM. Computations took less than 10 seconds for $N=512$ and less than 80 seconds for $N=2048$.

\subsection{Discussion}

A visible feature of Figures \ref{fig_recon} and \ref{fig_recon2}
is that the reconstructions do not change qualitatively by increasing the discretization
parameter $n$. This result is in line with Theorem \ref{pos_result}, i.e.,
if one fixes $\ep>0$ and takes $n$ to infinity then the minimizers converge
to Ambrosio--Tortorelli minimizers.

It is also evident that the parameter $\ep$ controls how sharp reconstructions
one can obtain. In Figure \ref{fig_recon2} this is visible with the second peak.
Namely, with the value $\ep=0.02$ this peak is smoothened whereas with the other values
the reconstruction becomes sharp. 

The convergence of the algorithm was satisfactory especially for $u$.
In most of the runs the value of $u$ was achieved very accurately with less than
10 iteration steps. However, the function $v$ still evolved slowly after this and 
a satisfactory estimate was obtained with 50 iteration steps where also each run was
stopped. The authors expect that this slowness can be overcome by more sophisticated minimization algorithm
in the step (3) of the algorithm. 

\appendix
\section{Technical lemmata}

\subsection{Properties of domain $X_\ep$}
\label{sec:choiceofa}

In the definition of domain $X_\ep$ in equation \eqref{Xep_domain} one 
restricts the values of function $v$ in the pair $(u,v)\in H^1(\T)\times H^1(\T)$
to the interval $[0,1+30\ep]$. Let us now discuss some properties related
to this choice. Define function
$g_\ep(t) = -\log(\ep^2+t^2)+\frac{1}{4\ep} (1-t)^2$
for $t\in\R$ where $\ep>0$ is fixed.
\begin{lemma}
  \label{aux_function}
  Assume that $0<\ep<\frac 1 8$.
  The function $g_\ep$ has a unique minimizer $t_\ep$ 
  which satisfies $1 \leq t_\ep \leq 1+ 30\ep$.
  Furthermore, the inequality $g_\ep(t) \leq g_\ep(s)$  
  holds when $s$ and $t$ satisfy one of the following conditions:
  \begin{itemize}
    \item[(i)] $1\leq t \leq 1+ 30\ep$ and $s>1+30\ep$,
    \item[(ii)]  $t\in [0,1]$ and $s\leq -1$ or
    \item[(iii)] $t\in [0,1]$ and $s=-t$.
  \end{itemize}
\end{lemma}

\begin{proof}
Clearly, one has $g_\ep(t)\leq g_\ep(-t)$ for any $t\in \R$. This proves claim $(iii)$ and 
since $\lim_{t\to \infty}g_\ep(t)=\infty$ this also shows
that the global minimizer has to be located in $\R_+$.

 The derivative $Dg_\ep$ has the form
  $Dg_\ep(t) = \frac{-2 t}{\ep^2+t^2} + \frac{1}{2\ep}(t-1)$
  for $t\in\R$. 
  The first term is negative everywhere in $\R_+$.
  Since the second term increases linearly and is positive for $t>1$, the zeros of $Dg_\ep$ on $\R_+$
  have to 
  be greater than $1$. 
  Also since $\lim_{t\to\infty} Dg_\ep(t) = \infty$ and 
the first term is strictly decreasing 
 for $t>1$ the function $D_\ep$ has a unique zero $t_\ep$ for $t>1$. This yields the existence of a unique
 minimizer for $g_\ep$.
Furthermore, claim $(ii)$ can be easily deduced since $Dg_\ep(t)<0$ for $t\leq 1$ when $\ep<1/8$.

Let us now show an upper bound to $t_\ep$. 
  Apply inequality 
  $\frac{t}{\ep^2+t^2} \leq \frac{1}{t}$
  for $t>1$ to obtain a lower bound to function $Dg_\ep$. By solving equation
  $Dg_\ep(t) \geq -\frac{1}{t_+}+\frac{1}{2\ep}(t_+-1)=0$
  for $t_+>1$ one obtains a bound $t_\ep \leq t_+$.
  A short computation yields
  $t_+ = 1 + \frac 12(\sqrt{1+8\ep}-1) \leq 1 + 2\ep$.
  
  Finally, let us study the claim $(i)$. 
  From above it is evident that there exists a unique point $s_\ep>t_\ep$ such that
  $g_\ep(s_\ep) = g_\ep(1)$. In the following we show that $s_\ep<1+30\ep$.  
  Denote $h_\ep(t) = g_\ep(t) - g_\ep(1)$. Then we have that
  \begin{equation}
    h_\ep(t) = -\log \frac{\ep^2+t^2}{\ep^2+1} + \frac{1}{4\ep}(1-t)^2
    \geq 1-\frac{\ep^2+t^2}{\ep^2+1} + \frac{1}{4\ep}(1-t)^2
    \label{hp_aux_ineq}
  \end{equation}
  for any $t\geq1$ where we have used inequalities $-\log x \geq -x + 1$ for $x\geq 0$.
  The quadratic function on the right-hand side has a zero in $t_1=1$ and with a detailed calculation one can show that the second zero 
  satisfies $t_2<1+30\ep$ for $\ep<\frac 18$. 
\end{proof}

\subsection{Auxiliary bounds}

Here we show auxiliary technical lemmata.
Define
\[
  G_{\ep,n} (v_n,b) = \int_\T \left(-N \log(\ep^2+(Q_nv_n)^2) + \frac{b}{4\ep} (1-v_n)^2\right) dt
\]
where $v_n\in PL(n)$, $b,\ep>0$, $\alpha\in\R$ and $n\in\N$.

\begin{lemma}
\label{infFep}
For any $0<\ep<\frac 18$, $n\in\N$ and $b\geq 0$ there are constants $C$ and $C(b)$ such that
\[
  -C(b)\ep N^{2} \leq
  \inf_{v\in PL(n)} G_{\ep,n}(v,b)
  \leq -C(\sqrt{\ep}N-1).
\]
\end{lemma}

\begin{proof}
The upper bound for the infimum follows by setting $v\equiv 1+\sqrt{\ep}$
and using inequality $\log(1+x) \geq \frac 12 x$
for small $x>0$.
For the lower bound first notice that
\[
   -\log(\ep^2+(Q_nv)^2) \geq -2 \log \sqrt{\ep^2+(Q_nv)^2}
   \geq -2\log(\ep+\abs{Q_nv}) \geq -2(\ep+\abs{Q_nv}-1).
\]
Since $\int_\T |Q_n v|dx \leq \int_\T |v| dx$ it holds also that
\begin{equation}
  \label{aux_infFep_ineq}
  \int_\T -N \log(\ep^2+(Q_nv)^2) dx
  \geq \int_\T-2 N (\ep+\abs{v}-1) dx.
\end{equation}
Now denote $h_\ep(t) = -2 N (\ep+\abs{t}-1) + \frac{b}{4\ep} (1-t)^2$ for any $t\in\R$.
Clearly we have $h_\ep(-t) \geq h_\ep(t)$ for $t\geq 0$.
For positive values of $t$ function $h_\ep$ is quadratic function
$h_\ep(t) = -2N \ep -2 N(t-1) + \frac{b}{4\ep} (t-1)^2$
with respect to variable $t-1\geq -1$. The minimum of this function is obtained 
when $t-1=\frac{4}{b}N\ep$ and thus
$h_\ep(t) \geq -2N \ep - \frac 4 b N^{2} \ep$.
It is now easy to verify that
\[
   \int_\T-2 N (\ep+\abs{v}-1) + \frac{b}{4\ep} (1-v)^2 dt
   \geq -\frac 4b N^{2}\ep -2N \ep 
   \geq -C(b)N^{2}\ep.
\]
Together with the inequality \eqref{aux_infFep_ineq} this yields the claim.
\end{proof}

\begin{lemma}
\label{logconv}
Assume that a sequence $v_\ep \in H^1(\T ; [0,1+C\ep])$ satisfies
$\int_\T (1-v_\ep)^2dt \leq C'\ep$
for some constants $C,C'>0$. Then it follows that
$\lim_{\ep\to 0} \pone(v_\ep) = 0$.
\end{lemma}

\begin{proof}
Let us denote $E_\ep =\{ t\in \T \; | \; v_\ep(t) \leq \frac 12\}$ for $\ep>0$.
The Lebesgue measure of $E_\ep$ is bounded by $|E_\ep| \leq C \ep$
and thus
$\int_{E_\ep} |\log(\ep^2+v_\ep^2) | dt \leq 2C \ep \log \ep$
which converges to zero as $\ep\to 0$.
Denote $\tilde v_\ep = \max(v_\ep,\frac 12)$. Clearly also $v_\ep \to \onef$ in $L^2(\T)$ and hence
by the Lebesgue dominated convergence theorem
$\lim_{\ep \to 0} \pone(v_\ep)
  \leq \lim_{\ep \to 0} \left(\pone(\tilde v_\ep)+ 2C\ep \log \ep\right) = 0.$
This proves the statement.
\end{proof}

The following lemma is proved in \cite{Saksman} in more detail.

\begin{lemma}
  \label{dual_bound}
  For any $0<s<\frac 12$, $u\in L^1(a,b)\cap H^{-s}(a,b)$ with $a,b\in\R$ such that $b>a$
  we have
  $\abs{\int_a^b u dt} \leq C |b-a|^{\frac 12 - s}\norm{u}_{H^{-s}(a,b)}$.
\end{lemma}

\begin{proof}
  By \cite{LiMa} the dual space of $H^{-s}(a,b)$ is
  $\widetilde{H}^s(a,b) = \{ f \in H^s(\R) \; | \; {\rm supp}(f) \subset [a,b]\}$  
  with norm $\norm{f}_{\widetilde{H}^s(a,b)} = \norm{f}_{H^s(\R)}$.
  Furthermore, the mapping 
  $T:f \mapsto f \onef_{(a,b)}$
  is continuous in $H^t(\R)$ for any $-1/2<t<1/2$. 
  In particular, we have that the $\onef_{(a,b)}\in \widetilde{H}^s(a,b)$.
  Without losing any generality we assume $a = -b$.
  The Fourier transform of $\onef_{(-b,b)}$ satisfies 
  $\widehat{\onef_{(-b,b)}}(\xi) = C \frac{\sin b \xi}{\xi}$ and thus
  \begin{eqnarray*}
    \bigg|\int_{-b}^b u dt\bigg| & \leq & \norm{\onef_{(-b,b)}}_{\widetilde{H}^s(-b,b)}\norm{u}_{H^{-s}(-b,b)} \\
    & = & C \left(\int_{-\infty}^\infty \abs{\frac{\sin b \xi}{\xi}}^2 (1 + \abs{\xi})^{2s}d\xi\right)^{1/2} \norm{u}_{H^{-s}(-b,b)}  \\
    & \leq & C' b^{\frac 12-s} \norm{u}_{H^{-s}(-b,b)} 
  \end{eqnarray*}
  for some constant $C'>0$.
\end{proof}

\newpage
\thispagestyle{empty}

\renewcommand{\thesection}{7}

\end{document}